\documentclass[12pt]{amsart}
\usepackage[numbers, square]{natbib}
\usepackage{amssymb}
\usepackage{amsmath}
\usepackage{amsfonts}
\usepackage{geometry}
\usepackage{graphicx}
\usepackage{amsthm}
\usepackage{hyperref}
\usepackage[dvipsnames]{xcolor}

\providecommand{\U}[1]{\protect\rule{.1in}{.1in}}
\theoremstyle{plain}

\newtheorem{lemma}{Lemma}

\newtheorem{proposition}{Proposition}
\newtheorem{remark}{Remark}

\newtheorem{theorem}{Theorem}
\numberwithin{equation}{section}

\DeclareMathOperator{\supp}{supp}

\begin{document}
\title[Doubling inequalities and critical sets of Dirichlet eigenfunctions]{Doubling inequalities and critical sets of Dirichlet eigenfunctions}
\author{ Jiuyi Zhu}
\address{
 Department of Mathematics\\
Louisiana State University\\
Baton Rouge, LA 70803, USA\\
Emails: zhu@math.lsu.edu}
\thanks{Zhu is supported in part by NSF grant OIA-1832961. }
\date{}
 \dedicatory{}

\begin{abstract}
We study the sharp doubling inequalities for the gradients and  upper bounds for the critical sets of Dirichlet eigenfunctions on the boundary and in the interior of compact Riemannian manifolds. Most efforts are devoted to obtaining the sharp doubling inequalities for the gradients.
 New technique is developed to overcome the difficulties on the unavailability of the double manifold in obtaining doubling inequalities in smooth manifolds. The sharp upper bounds of critical sets in analytic Riemannian manifolds are consequences of the doubling inequalities.
\end{abstract}

\maketitle
\section{Introduction}
In this paper, we mainly study the doubling inequalities for the gradients 
of Dirichlet eigenfunctions
\begin{equation}
\left\{
\begin{array}{lll}
-\triangle_g e_\lambda(x)=\lambda e_\lambda(x),
 \quad &x\in\mathcal{M},
\medskip \\
 e_\lambda(x)=0,
\quad &x\in
\partial\mathcal{M},
\end{array}
\right. \label{diri}
\end{equation}
where $(\mathcal{M}, g)$ is a smooth and compact $n$-dimensional Riemannian manifold with boundary ($n\geq 2$). With the sharp doubling inequalities of the gradients on the boundary and in the interior of the manifold, we are able to obtain the sharp upper bounds of critical sets of Dirichlet eigenfunctions by standard techniques.
Doubling inequalities are inequalities with norms in a ball controlling the norms in a double ball. They are effective tools to control the local growth of the functions. Thus, these inequalities are important in the quantitative characterization of the strong unique continuation property and in the study of the measure of zero-level sets of given functions.

 Let us briefly review the relations of doubling inequalities and the measure of nodal sets (zero-level sets) of eigenfunctions. The doubling inequalities have played  important roles in the measure of nodal sets of eigenfunctions
\begin{equation}
-\triangle_g \phi_\lambda=\lambda \phi_\lambda
\label{class}\end{equation}
 on compact manifolds $(\mathcal{M}, g)$
without boundary. The celebrated problem about nodal sets centers
around the famous Yau's conjecture for smooth manifolds.
Yau \cite{Y} conjectured that  the upper and lower bound of nodal sets of
eigenfunctions in (\ref{class}) are controlled by
\begin{equation}c\sqrt{\lambda}\leq H^{n-1}(\{x\in\mathcal{M}|\phi_\lambda(x)=0\})\leq C\sqrt{\lambda}\label{yau}
\end{equation} where $C, c$
depend only on the manifold $\mathcal{M}$. The conjecture was shown
to be true for real analytic manifolds by Donnelly-Fefferman in
\cite{DF}, \cite{DF1}. Lin \cite{Lin} also proved the upper bound for the
analytic case using a different approach. The following doubling inequality 
\begin{align}
\|\phi_\lambda \|_{L^\infty(\mathbb B_{2r}(x))}\leq e^{{C} \sqrt{\lambda}} \|\phi_\lambda \|_{L^\infty(\mathbb B_{r}(x))}
\label{kaodouble}
\end{align}
for any $B_{r}(x)\subset \mathcal{M}$ is essential in the derivations of bounds (\ref{yau}). Note that (\ref{kaodouble}) is obtained in the smooth manifold.
For the conjecture (\ref{yau}) on the measure of nodal sets on smooth manifolds, there are important breakthrough made by Logunov and Malinnikova \cite{LM},  \cite{Lo1} and \cite{Lo2} in recent years.

It is well known that eigenfunctions $e_\lambda$ change signs as eigenvalues $\lambda$ increase. The nodal sets of eigenfunctions may intersect other interior nodal sets or may intersect the boundary. Critical sets of eigenfunctions are sets where the gradients of eigenfunctions $\nabla e_\lambda$ vanish. The critical sets appear at the intersections of nodal sets or inside of nodal domains. We aim to study the measure of critical sets on the boundary $\{x\in \partial\mathcal{M}| |\nabla e_\lambda|=0\}$, where $\nabla e_\lambda$ is the full derivative in $ {\mathcal{M}}$. Since $e_\lambda=0$ on boundary, those critical sets on the boundary can be regarded as  singular sets, that is, $\{x\in \partial\mathcal{M}| e_\lambda=|\nabla e_\lambda|=0\}$. By Hopf lemma, we can  check that these critical sets are  the intersection of nodal sets in the interior manifold with the boundary. More intuitively, these critical sets are where interior nodal sets touch the boundary. We can also consider these sets as the intersection of different nodal sets if we can extend the Dirichlet eigenfunction across the boundary.  Since the boundary $\partial\mathcal{M}$ is $(n-1)$ dimensional, the boundary critical sets $\{x\in \partial\mathcal{M}| |\nabla e_\lambda|=0\}$ are with dimension no more than $(n-2)$. To obtain the measure of boundary critical sets, we study the measure of zero sets of the normal derivatives $\frac{\partial e_\lambda}{\partial \nu}$ on the boundary, where $\nu$ is a unit outer normal on the boundary. Based on the relations of doubling inequalities and zero-level sets, in order to obtain the measure of boundary critical sets, we need to obtain the doubling inequality for normal derivatives on the boundary for smooth manifolds.
\begin{theorem}
Let $e_\lambda$ be the Dirichlet eigenfunctions in (\ref{diri}). There exist positive constants $C$ and $r_0$ depending only on the smooth manifold $\mathcal{M}$ such that
\begin{equation}
\|\frac{\partial e_\lambda}{\partial \nu}\|_{L^\infty(\mathbb B_{2r}(x))}\leq e^{C\sqrt{\lambda}}\|\frac{\partial e_\lambda}{\partial \nu}\|_{L^\infty(\mathbb B_{r}(x))}
\label{LLLkao}
\end{equation}
for any $0<r<r_0$ and any $\mathbb B_{2r}(x)\subset\partial\mathcal{M}$. 
\label{th0}
\end{theorem}

Together with zeros counting results for complex analytic functions, 
we can show the following sharp upper bounds for the boundary critical sets.
\begin{theorem}  \label{th1}
Let $\mathcal{M}$ be a real analytic manifold
with boundary and $e_\lambda$ be the Dirichlet eigenfunctions in (\ref{diri}). There exists a positive constant $C$ depending on $\mathcal{M}$
such that
\begin{align}  H^{n-2}(\{x\in \partial \mathcal{M}| \ |\nabla e_\lambda|=0  \} )\leq C\sqrt{\lambda}.
 \end{align}
\end{theorem}
For planar analytic domains, such upper bounds for boundary critical sets of  Dirichlet eigenfunctions was derived by Toth and Zelditch \cite{TZ} among other interesting results using a different method. Theorem \ref{th1} improves such upper bounds to general dimensions in the analytic manifold. See also \cite{TZ1} for their related interesting work on nodal intersections with a certain analytic hypersurface in the analytic manifold without boundary.

Next we focus on the study of interior critical sets and interior doubling inequalities.
For elliptic equations without zero order term, the critical sets of the solutions are at most $(n-2)$-dimensional. Exponential upper bounds have been established for Hausdorff measure of critical sets of solutions in smooth manifolds by Naber and Valtorta \cite{NV}, see also other work on the study of critical sets or singular sets of elliptic equations, \cite{CNV}, \cite{Han}, \cite{Han1}, \cite{HHL}, \cite{HHON}, just name a few.  For the Laplace eigenfunctions, the dimension of critical sets of eigenfunction can be varied. It is shown by Jackobson and Nadirashvili \cite{JN} that there exists a Riemannian metric on two dimensional torus $\mathbb T^2$ such that
the number of critical points is uniformly bounded for some sequence of eigenfunctions, which disproves the Yau's conjecture \cite{Y} on the increase of the number of critical points. 
Recently, Buhovsky, Logunov and Sodin \cite{BLS} constructed a Riemannian metric on the two dimensional torus $\mathbb T^2$, such that for infinitely many eigenvalues, each corresponding eigenfunction has infinitely many isolated critical points.

There exist examples of Dirichlet eigenfunctions which have $(n-1)$-dimensional critical sets.  Let $(\mathcal{N}, g)$ be a $(n-1)$-dimensional compact manifold and $\mathcal{M}=[0, 1]\times \mathcal{N}$. The compact manifold $\mathcal{M}$ is equipped with the product metric. Then the boundary $ \mathcal{M}$ is given by $\partial\mathcal{M}=\{0, 1\}\times \mathcal{N}$. We can check that $e_k(x)=\sin(2\pi kx)$ are Dirichlet eigenfunctions with eigenvalues $\lambda= 4\pi^2 k^2$ on $\mathcal{M}$. We can also see that the critical set is $(n-1)$-dimensional. Furthermore, $c\sqrt{\lambda}\leq H^{n-1}(\{x\in \mathcal{M}||\nabla e_k|=0  \} )\leq C\sqrt{\lambda}$, where $c, C$ depend on the manifold $\mathcal{M}$.  There are also some surfaces of revolution with $(n-1)$-dimensional critical sets, see also \cite{Z}. Since the critical sets of Dirichlet eigenfunction are at most $(n-1)$-dimensional, we can bound the Hausdorff measure of critical sets as we do for zero-level sets.
We first need to obtain the doubling inequality for the gradients of Dirichlet eigenfunctions near the boundary. For Dirichlet eigenfunctions, even if we can do an odd extension to have a double manifold, we can not directly show the doubling inequalities for $\nabla e_\lambda$, because the metric for the double manifold is only Lipschitz. To avoid the use of the double manifold, we develop the quantitative Carleman estimates in the half balls,  and then show new propagation arguments to obtain the doubling inequalities.
\begin{theorem}
Let $e_\lambda$ be the Dirichlet eigenfunction in (\ref{diri}). There exist positive constants $C$ and $r_0$ depending only on the smooth manifold $\mathcal{M}$ such that
\begin{equation}
\|\nabla e_\lambda\|_{L^\infty(\mathbb B^+_{2r}(x_0))}\leq e^{C\sqrt{\lambda}}\|\nabla e_\lambda\|_{L^\infty(\mathbb B^+_{r}(x_0))}
\label{LLL}
\end{equation}
for any $0<r<r_0$ and any $\mathbb B^+_{2r}(x_0)\subset\mathcal{M}$ and $x_0\in \partial\mathcal{M}$.
\label{th3}
\end{theorem}

Thanks to the doubling inequality in half balls (\ref{LLL}), and the interior doubling inequality in Lemma \ref{submit},  and zeros counting results, the upper bounds of interior critical sets are derived.
\begin{theorem}
Let $\mathcal{M}$ be a real analytic manifold
with boundary and $e_\lambda$ be the Dirichlet eigenfunctions in (\ref{diri}). There exists a positive constant $C$ depending on $\mathcal{M}$
such that
\begin{align}  H^{n-1}(\{x\in \mathcal{M}||\nabla e_\lambda|=0  \} )\leq C\sqrt{\lambda}.
 \end{align}
   \label{th2}
\end{theorem}
For real analytic manifolds without boundary, the same strategy was used in \cite{B} to obtain the upper bounds of critical sets of eigenfunctions in (\ref{class}). In our proof of Theorem \ref{th2}, we need to overcome the difficulties on the presence of the boundary. Specifically, we have to deal with the critical sets in the neighborhood of the boundary and in the interior of the manifold. Furthermore, new challenges arise when dealing with the doubling inequalities of the gradients of Dirichlet eigenfunctions.

The following remark provides new insights for the techniques used in the proof of Theorem \ref{th3}.
\begin{remark}
The method in the proof of Theorem \ref{th3} provides new ways to prove the doubling inequality for smooth manifolds if the double manifold is not available. It has interesting applications to the study of nodal sets or other quantitative properties of eigenfunctions of general eigenvalue problems. For example,
the Steklov problem is often
considered in this more general form with physical background
\begin{equation}
\left \{ \begin{array}{lll}
\triangle_g e_\lambda=0 \quad &\mbox{in} \ {\mathcal{M}},  \medskip\\
\frac{\partial e_\lambda}{\partial \nu}=\lambda \varrho(x) e_\lambda\quad &\mbox{on} \ {\partial\mathcal{M}}.
\end{array}
\right.
\label{steklov}
\end{equation}
 The weight function $\varrho(x)$ represents the mass density along the boundary $\partial\mathcal{M}$. Another example is a general type of Robin eigenvalue problem with weight functions on the boundary
\begin{equation}
\left \{ \begin{array}{lll}
-\triangle_g e_\lambda=\lambda e_\lambda \quad &\mbox{in} \ {\mathcal{M}},  \medskip\\
\frac{\partial e_\lambda}{\partial \nu}=\alpha V(x) e_\lambda\quad &\mbox{on} \ {\partial\mathcal{M}}.
\end{array}
\right.
\label{Robin}
\end{equation}

 The double manifold does not seem to be available because of the presence the function $\varrho(x)$ or $V(x)$. The technique in Theorem \ref{th3} provides the way to prove the doubling inequality for the half balls centered at the boundary for $e_\lambda$ in (\ref{steklov}) or (\ref{Robin}). We can further prove the upper bounds of interior nodal sets of $e_\lambda$ in the manifold. See also other work e.g. \cite{BL}, \cite{Z1}, \cite{Zh1}, \cite{Zh2}, \cite{Zh3}, \cite{LZ} on the upper bounds of interior nodal sets or boundary nodal sets for eigenfunctions as (\ref{steklov}) or (\ref{Robin}).
\end{remark}

 The outline of the paper is as follows. Section 2 is devoted to the proof of Theorem \ref{th0} and Theorem \ref{th1}.
 In section 3, we first derive the systems of equations which gradients of eigenfunctions $\nabla e_\lambda$ satisfy. Then we state the quantitative Carleman estimates on the half balls. As applications of Carleman estimates, we prove three half-ball inequalities and Theorem \ref{th3} in section 4. New propagation arguments are applied to bypass the use of the double manifold. In section 5, we derive the proof of Theorem \ref{th2}. The last section is for the proof of the the quantitative Carleman estimates in half balls.
 The letter $c$, $C$, $C_i$ denote generic
positive constants and do not depend on $\lambda$. They may vary in
different lines and sections. In the paper, since we study the asymptotic properties for eigenfunctions, the eigenvalue $\lambda$ is assumed to be sufficiently large.
\medskip

\noindent{\bf Acknowledgement.}  The author appreciates Professor Fang-hua Lin for helpful discussions on the topics. The author also thanks Professor  Steve Zelditch for bringing \cite{TZ1} to our attentions and for information on the possible extension of their results to critical sets.

\section{Boundary critical sets}
In this section, we first prove the sharp doubling inequality for the normal derivatives, then show the upper bounds of the measure of critical sets on the boundary for Dirichlet eigenfunctions.
 In the local coordinate charts, using the Einstein notations (the summation notation is understood), the equation is written as
\begin{align}
-\triangle_g e_\lambda= - g^{-\frac{1}{2}} \partial_i (g^{\frac{1}{2}} g^{ij}  \frac{\partial e_\lambda}{\partial x_j})=\lambda e_\lambda,
\end{align}
where $g^{ij}$ denotes the inverse of the metric $g_{ij}$ and  $g=\det(g_{ij})$. We will write $\triangle_g$ as $\triangle$ if the context is understood. We need to pay special attention to the equation around the boundary $\partial\mathcal{M}$.  For any point $p\in \partial\mathcal{M}$, the Fermi exponential map at $p$ which gives the Fermi coordinate systems, is defined on a half ball of $\mathbb R^n_+ \approx T_p(\partial \mathcal{M})$ centered at origin. Suppose $(x_1, \cdots, x_{n-1})$ is the geodesic normal coordinate of $\partial\mathcal{M}$ at $p$. Let $x_n =dist (x, \ \partial \mathcal{M})$. Note that $dist (x, \ \partial \mathcal{M})$ is smooth in a small open neighborhood of $p$ in ${ \mathcal{M}}$ if $\mathcal{M}$ is smooth. We can locally identify $\partial \mathcal{M}$ as $x_n=0$.
 By Fermi geodesic coordinates,
\begin{align} -g^{-\frac{1}{2}} \partial_i ( g^{\frac{1}{2}} g^{ij} \frac{\partial e_\lambda}{\partial x_j})=\lambda e_\lambda
\label{fermi}
\end{align}
with $g^{nn}=1$, $g^{in}=0$ and $g^{ij}(x', x_n)\not =0$ for $1\leq i, j\leq n-1$.

We will take advantage of the new quantitative propagation smallness results for the second order elliptic equations
in the half ball shown in \cite{Zh3}, that is, quantitative two half-ball and one lower dimensional ball inequality. Let us present the results in a general setting. Let $u$ be the solutions of
\begin{align}
- a_{ij} D_{ij}u +b_i(x) D_i u +c(x)u=0  \quad{in} \ \mathbb B^+_{1/2},
\label{general}
\end{align}
where $a_{ij}$ is $C^1$, $ b(y)$ and $ c(y)$ satisfy
\begin{equation}
\left\{ \begin{array}{lll}
\| {b}\|_{W^{1, \infty} (\mathbb B^+_{1/2})}\leq C(\tau_1+1), \medskip \\
\| {c}\|_{W^{1, \infty}(\mathbb B^+_{1/2})}\leq C(\tau_1^2 +\tau_2),
\end{array}
\right.
\label{aaa}
\end{equation}
and $\tau_1$ and $\tau_2$ are positive constants with possible large values.
The quantitative two half-ball and one lower dimensional ball inequality is stated as follows.
\begin{lemma}
Let $u\in C^\infty_0(\mathbb B^+_{1/2}  )$ be a solution of (\ref{general}). Denote the lower dimensional ball
$${\bf B}_{1/3} =\{ (x', \ 0)\in \mathbb R^n| x'\in \mathbb R^{n-1}, \ |x'|<\frac{1}{3}\}.$$
Assume that
\begin{equation}
\|u \|_{H^1({\bf B}_{1/3})}+\|\frac{\partial u}{\partial \nu}\|_{L^2({\bf B}_{1/3})}\leq \epsilon <<1
\end{equation}
and $\|u \|_{L^2(\mathbb B^+_{1/2})}\leq 1$. There exist positive constants $C$ and $\beta$ such that
\begin{equation}
\|u\|_{L^2(\frac{1}{256}\mathbb B^+_1)}\leq e^{C(\tau_1+\sqrt{\tau_2})} \epsilon^\beta.
\label{too}
\end{equation}
More precisely, we can show that there exists $0<\kappa<1$ such that
\begin{equation}
\|u\|_{L^2(\frac{1}{256}\mathbb B^+_1)}\leq e^{C(\tau_1+\sqrt{\tau_2})} \|u\|_{L^2(\mathbb B^+_{1/2})}^\kappa \big( \|u \|_{H^1({\bf B}_{1/3})}+\|\frac{\partial u}{\partial \nu}\|_{L^2({\bf B}_{1/3})}  \big)^{1-\kappa}.
\label{tool}
\end{equation}
\label{halfspace}
\end{lemma}
 These qualitative results were established in \cite{Lin} and \cite{ARRV}. We obtained the quantitative results with the consideration of the quantitative behavior of $\tau_1$ and $\tau_2$  by some global Carleman estimates. The estimates (\ref{tool}) have already shown their important applications on the measure of boundary nodal sets of Neumann eigenfunctions in \cite{Zh3}.

 For the Dirichlet eigenfunctions, the following sharp doubling inequalities hold
\begin{align}
\|e_\lambda \|_{L^2(\mathbb B^+_{2r}(x))}\leq e^{{C} \sqrt{\lambda}} \|e_\lambda \|_{L^2(\mathbb B^+_{r}(x))}
\label{drdouble}
\end{align}
 in half balls for any $0< r\leq r_0$, where $r_0$ depend only on $\mathcal{M}$ and $x\in \partial\mathcal{M}$. The estimate (\ref{drdouble}) was actually established in \cite{DF}, even if it was not stated for the half balls, because the doubling inequality in balls in the double manifold was sufficient for their results.  Note that the sharp doubling inequalities (\ref{drdouble}) can also be proved by the technique developed in section 4 without using the technique of the double manifold.

Together with Lemma \ref{halfspace} and doubling inequalities (\ref{drdouble}) in half balls, we can readily derive the doubling inequality for the gradients of Dirichlet eigenfunctions on the bounday of the manifold.

\begin{proof}[Proof of Theorem \ref{th0}]
Because of the Fermi coordinates, we consider the Dirichlet eigenfunctions $ e_\lambda$ in (\ref{fermi}) near the boundary. We may argue on scale of order one
and  normalize $e_\lambda$ as
\begin{align} \|e_\lambda \|_{L^2(\mathbb B^+_{1/2})}=1.
\label{normalize} \end{align}

 Note that $e_\lambda=\nabla_t e_\lambda=0$ on the boundary $\{x_n=0\}$, where $\nabla_t$ is the derivative in the tangential direction. We can write the equation (\ref{fermi}) in the form of (\ref{general}). In term of the assumptions (\ref{aaa}), we can choose $\tau_1$ as some fixed constant depending on $\mathcal{M}$ and $\tau_2=\lambda$. Hence the quantitative three-ball inequality on the half balls (\ref{tool}) holds. We may normalize it  as
\begin{equation}
\|e_\lambda \|_{L^2(\mathbb B^+_{1/512})}\leq e^{ C\sqrt{\lambda}} \|e_\lambda \|_{L^2(\mathbb B^+_{1/4})}^\kappa  \|\frac{\partial e_\lambda}{\partial x_n} \|^{1-\kappa}_{L^2({\bf B}_{1/6})}  .
\label{tooll}
\end{equation}
By finitely many iterations of the doubling inequality in the half balls (\ref{drdouble}) and (\ref{normalize}), we get
\begin{align}
\|\frac{\partial e_\lambda}{\partial x_n}  \|_{L^2({\bf B}_{1/6})} \geq e^{ -C\sqrt{\lambda}}.
\label{haoma}
\end{align}

Let $\eta$ be a cut-off function  such that $\eta(y)=1$ for $|y|\leq \frac{1}{4}$ and vanishes for $|y|\geq \frac{1}{3}$.  By the Hardy trace inequality and elliptic estimates, it follows that
\begin{align}
\|\frac{\partial e_\lambda}{\partial x_n} \|_{L^2( {\bf B}_{1/4})} &\leq \| \eta \frac{\partial e_\lambda}{\partial x_n}\|_{L^2( {\bf B}_{1/4})}  \leq \|\nabla( \eta \frac{\partial e_\lambda}{\partial x_n} )\|_{L^2( \mathbb R^{n}_+)} \nonumber \\
&\leq C \|\nabla^2 e_\lambda\|_{L^2( {\mathbb B}^+_{1/3})}+ C \|\frac{\partial e_\lambda}{\partial x_n} \|_{L^2( {\mathbb B}^+_{1/3})}\nonumber \\
&\leq C{\lambda} \| e_\lambda \|_{L^2( {\mathbb B}^+_{1/2})}\nonumber \\
&\leq C{\lambda}.
\label{lastt}
\end{align}
Combining the established estimates (\ref{haoma}) and (\ref{lastt}), we have
\begin{equation}
\|\frac{\partial e_\lambda}{\partial x_n} \|_{L^2( {\bf B}_{1/4})}\leq e^{ C\sqrt{\lambda}} \|\frac{\partial e_\lambda}{\partial x_n}\|_{L^2( {\bf B}_{1/6})}.
\label{night}
\end{equation}
From elliptic estimates and the way to derive (\ref{haoma}), the estimate (\ref{night}) holds in $L^\infty$ norm,
\begin{equation}
\|\frac{\partial e_\lambda}{\partial x_n} \|_{L^\infty( {\bf B}_{1/5})}\leq e^{ C\sqrt{\lambda}} \|\frac{\partial e_\lambda}{\partial x_n}\|_{L^\infty( {\bf B}_{1/6})}.
\label{nightba}
\end{equation}
 By rescaling and iteration, we arrive at
\begin{align}
\|\frac{\partial e_\lambda}{\partial x_n }\|_{L^\infty( {\mathbb B}_{2r}(x_0))}  \leq e^{ C\sqrt{\lambda}} \|\frac{\partial e_\lambda}{\partial x_n }\|_{L^\infty( {\mathbb B}_{r}(x_0))}
\label{double}
\end{align}
for any $x_0\in \partial\mathcal{M}$ and ${\mathbb B}_{2r}(x_0)\subset \partial\mathcal{M}$, and $r<r_0$ for some $r_0$ depending only on $\partial\mathcal{M}$.  This completes the proof of Theorem \ref{th0}.
\end{proof}

To measure the number of zeros, we need a lemma concerning the growth of a complex analytic function. See e.g. Lemma 2.3.2 in \cite{HL}.
\begin{lemma}
Suppose $f: \mathcal{B}_1(0)\subset \mathbb{C}\to \mathbb{C}$ is an
analytic function satisfying
$$ f(0)=1\quad \mbox{and} \quad \sup_{ \mathcal{B}_1(0)}|f|\leq 2^{\hat{N}}$$
for some positive constant $\hat{N}$. Then for any $r\in (0, 1)$, there
holds
$$\sharp\{z\in\mathcal{B}_r(0): f(z)=0\}\leq c\hat{N}       $$
where $c$ depends on $r$. Especially, for $r=\frac{1}{2}$, there
holds
$$\sharp\{z\in \mathcal{B}_{1/2}(0): f(z)=0\}\leq \hat{N}.        $$
\label{wwhy}
\end{lemma}

With the doubling inequality (\ref{LLLkao}) and Lemma \ref{wwhy} on hand, it is kind of standard to provide the proof of upper bounds for the boundary critical sets of Dirichlet eigenfunctions. See pioneer work on the measure of nodal sets by \cite{DF}, \cite{Lin} and other related work using this idea, e.g. \cite{Zh2}, \cite{Zh3}, \cite{LZ}
\begin{proof}[Proof of Theorem \ref{th1}] To get the hypoelliptic  estimates for Dirichlet eigenfunctions on the boundary, we perform a standard lifting argument. Let
\begin{align}
{w}(x, t)=e^{\sqrt{\lambda} t}e_\lambda(x).
\label{lift}
\end{align}
Then ${w}(x, t)$ satisfies the following equation
\begin{equation}
\left \{ \begin{array}{rll}
{\triangle}  {w}+\partial^2_t w=0 \quad &\mbox{in} \ {\mathcal{M}}\times (-\infty, \ \infty),  \medskip\\
w=0 \quad &\mbox{on} \ {\partial\mathcal{M}}\times (-\infty, \ \infty).
\end{array}
\right.
\label{error}
\end{equation}
 By straightening the boundary $\partial\mathcal{M}$ locally, rescaling and translation. we may assume that $(p, 0, t)\in \big(\partial \mathbb B^+_{1/16}\cap \{x_n=0\}\big)\times (-\frac{1}{16}, \frac{1}{16}) $ with $p\in \mathbb R^{n-1}.$  From elliptic estimates in Lemma 2.3 in \cite{MN},  we obtain that
\begin{align}
| \frac{\nabla^{\bar\alpha} {w} (p, 0, 0)}{ \alpha !(\alpha_n+1)} |\leq C \hat{C}^{k+1}  \|{w}   \|_{L^\infty \big(\mathbb B^+_{1/8}\times(-\frac{1}{8}, \frac{1}{8})\big)},
\end{align}
where
$ \bar \alpha=(\alpha_1, \cdots, \alpha_{n-1}, \alpha_n+1, 0)$ and $|\alpha|=|(\alpha_1, \cdots, \alpha_{n-1}, \alpha_n)|=k$, and $C$, $\hat{C}>1$ depends on $\mathcal{M}$.
By the definition of ${w}$ and (\ref{haoma}), we have that
 \begin{align}
| \frac{\nabla^{\alpha} \frac{\partial e_\lambda}{\partial x_n} (p, 0)}{ \alpha !} |&\leq C (\alpha_n+1) \hat{C}^{k}  \|w \|_{L^\infty \big(\mathbb B^+_{1/8}\times(-\frac{1}{8}, \frac{1}{8})\big)} \nonumber \\
& \leq C \hat{C}^{k} e^{C\sqrt{\lambda}} \| e_\lambda  \|_{L^\infty (\mathbb B^+_{1/8})} \nonumber \\
&\leq C \hat{C}^{k} e^{C\sqrt{\lambda}} \|  {\frac{\partial e_\lambda}{\partial x_n}}  \|_{L^\infty ({\bf B}_{1/4})}.
\end{align}
Then $ \frac{\partial e_\lambda}{\partial x_n}(p, 0)$ is real analytic for any $(p, 0)\in \partial \mathbb B^+_{1/16}\cap\{x_n=0\}$. We may consider $p$ as the origin in $\mathbb R^{n-1}$. By summing up a geometric series, we derive a holomorphic extension of $\frac{\partial e_\lambda}{\partial x_n}$ with
\begin{equation}
\sup_{ |z|\leq  \frac{1}{2\hat{C}} }| \frac{\partial e_\lambda}{\partial x_n}(z)|\leq e^{C\sqrt{\lambda}}  \|  {\frac{\partial e_\lambda}{\partial x_n}}  \|_{L^\infty ({\bf B}_{1/4})},
\label{zeld}
\end{equation}
where $\frac{1}{2\hat{C}}<\frac{1}{8}$ and $z\in \mathbb C^{n-1}$.
Note that ${\bf B}_r$ is denoted as the ball in $\mathbb R^{n-1}$ with radius $r$.
Taking the boundary doubling inequality (\ref{night}) and (\ref{zeld}) into consideration, by finitely many steps of iterations,
we conclude that
\begin{align}
\sup_{ |z|\leq  \frac{1}{2\hat{C}}}| \frac{\partial e_\lambda}{\partial x_n} (z)|\leq e^{C\sqrt{\lambda}}\sup_{ x\in {\bf B}_{\frac{1}{4\hat{C}}}}|  \frac{\partial e_\lambda}{\partial x_n} (x)|.
\end{align}
By rescaling arguments, we further derive that
\begin{align}
\sup_{ |z|\leq 2 r}| \frac{\partial e_\lambda}{\partial x_n} (z)|\leq e^{C\sqrt{\lambda}}\sup_{ x\in  {\bf B}_r}| \frac{\partial e_\lambda}{\partial x_n} (x)|,
\label{eventu}
\end{align}
where $0<r<\hat{r}_0$ and $\hat{r}_0$, $C$ depend on $\mathcal{M}$.

For ease of notations, let $v=\frac{\partial e_\lambda}{\partial x_n}$.
Thanks to the doubling inequality (\ref{eventu}) and the growth control lemma for zeros, i.e. Lemma \ref{wwhy}, we can  give the proof of Theorem \ref{th1}.
We may also argue on scales of one as well. Let $p\in {\bf B}_{1/4}\subset \mathbb R^{n-1}$ be the point where the supremum of $| v |$ is achieved. After rescaling, we can assume that $|v(p)|=1$.
For each direction $\omega\in S^{n-2}$, we consider the function $$ v_\omega(z)= v(p+z\omega), \quad \quad z\in \mathcal{B}_1(0)\subset \mathbb{C}.$$ Denote $N(\omega)=\sharp\{z\in\mathcal{B}_{1/2}(0)\subset
\mathbb{C}| v_\omega(z)=0\}$. With aid of the doubling inequality (\ref{eventu}) and Lemma \ref{wwhy}, we can show that
\begin{align}
&\sharp\{ x\in {\bf B}_{1/2}(p)\subset \mathbb R^{n-1} | x-p \ \mbox{is parallel to} \ \omega \ \mbox{and} \  v(x)=0\} \nonumber \\
&\leq \sharp\{ z\in \mathcal{B}_{1/2}(0)\subset \mathbb{C}| v_\omega(z)=0\}  \nonumber \\
&=N(\omega)  \nonumber \\
&\leq C\sqrt{\lambda}.
\end{align}
 By the integral geometry estimates, we readily deduce that
\begin{align}
H^{n-2}( \{ x\in {\bf B}_{1/2}(p)|  \frac{\partial e_\lambda}{\partial x_n} (x)=0\})
&\leq \int_{S^{n-2}}N(\omega)\, d\omega \nonumber\\
 &\leq \int_{S^{n-2}} C\sqrt{\lambda}\, d\omega \nonumber\\
&\leq C\sqrt{\lambda}.
\end{align}
Thus, we obtain the upper bound of critical sets
\begin{align}
H^{n-2}( \{ x\in {\bf B}_{1/4}(0)| \  |\nabla e_\lambda (x)|=0\}) \leq C\sqrt{\lambda}.
\end{align}
By rescaling, it also implies that
\begin{align}
H^{n-2} (\{ \mathbb B_{r_0}(p)\subset \partial\mathcal{M} |\ |\nabla e_\lambda (x)|=0\}) \leq C\sqrt{\lambda}
\end{align}
for some $r_0$ depending only on $\mathcal{M}$ and for any $p\in \partial\mathcal{M}$.
 Since the boundary $\partial\mathcal{M}$ is compact, by finite number of coverings, we complete the proof of the theorem.
\end{proof}

\section{Quantitative Carleman estimates}
In this section, we first derive the systems of elliptic equations for the derivatives of Dirichlet eigenfunctions. Then we will establish the quantitative Carleman estimates for the elliptic systems.
To study the interior critical sets, we take derivative $\partial_k$ with respect to $x_k$ for each $k=1, \cdots n$ in each coordinate chart. Thus, we have
\begin{align}
 &-  \partial_k g^{-\frac{1}{2}}  \partial_i (g^{\frac{1}{2}} g^{ij}  \frac{\partial e_\lambda}{\partial x_j})
 -   g^{-\frac{1}{2}} \partial_i ( \partial_k ( g^{\frac{1}{2}}g^{ij} )\frac{\partial e_\lambda}{\partial x_j})
  -g^{-\frac{1}{2}} \partial_i (  (g^{\frac{1}{2}} g^{ij}  ) \partial_k\frac{\partial e_\lambda}{\partial x_j})
 =\lambda \partial_k e_\lambda.
 \label{deriv}
\end{align}
Let $U=\langle U_1, U_2\cdots, U_n \rangle$, where each $U_k=\frac{\partial e_\lambda}{\partial x_k}$ for $k=1, \cdots, n$. We rearrange (\ref{deriv}) as follows
\begin{align}
\triangle_g U_k +\langle B_k,  \ \nabla U\rangle +A_k \cdot U+\lambda U_k=0
\end{align}
where $A_k$ is a vector function and $ B_k$ is a matrix function depending only on the metric $g$ and its derivatives. Hence, we have a system of equation for
$U$,
\begin{align}
\triangle_g U +\langle B, \  \nabla U\rangle +A \cdot U+\lambda U=0,
\label{system1}
\end{align}
where $B$ a matrix consisting of $B_k$ and $A$ is a matrix consisting of $A_k $.

We also need to derive the boundary conditions for elliptic systems of the derivatives of Dirichlet eigenfunctions.
From the Dirichlet boundary conditions, in the local coordinate, we derive that
\begin{align}
\frac{\partial e_\lambda}{\partial x_1}=\cdots=\frac{\partial e_\lambda}{\partial x_{n-1}}=0 \quad \mbox{on} \ x_n=0.
\label{boundary23}
\end{align}
Thus, we know $U_1=\cdots =U_{n-1}=0$ on $x_n=0$. We also need to deduce another boundary condition for $U_{n}=\frac{\partial e_\lambda}{\partial x_{n}}$.
We write the equation (\ref{fermi}) in local Fermi coordinates as
\begin{align}
-g^{-\frac{1}{2}} \partial_i ( g^{\frac{1}{2}} g^{ij} ) \frac{\partial e_\lambda}{\partial x_j}- g^{ij}\frac{\partial^2 e_\lambda}{\partial x_{i}\partial x_j} =\lambda e_\lambda.
\end{align}
Since the equation is smooth up to the boundary, from the boundary conditions (\ref{boundary23}) and properties of Fermi coordinates, we obtain that
\begin{align}
-g^{-\frac{1}{2}} \frac{\partial}{\partial x_n} ( g^{\frac{1}{2}} g^{nn} ) \frac{\partial e_\lambda}{\partial x_n}- g^{nn}\frac{\partial^2 e_\lambda}{\partial x^2_{n}}   =\lambda e_\lambda \quad \mbox{on} \ x_n=0.
\end{align}
By the fact that $e_\lambda=0$ on $x_n=0$, we show that
\begin{align}
\frac{\partial U_n}{\partial x_n}=- \frac{\partial  \ln g^{\frac{1}{2}}}{\partial x_n} U_n.
\end{align}
Thus, to study the doubling inequalities and interior critical sets of Dirichlet eigenfunctions, it is reduced to investigate the elliptic systems with boundary conditions as follows
\begin{equation}
\left \{ \begin{array}{lll}\triangle_g U +\langle B, \  \nabla U\rangle +A \cdot U+\lambda U=0, \quad \quad &\mbox{in} \ \mathbb B^+_1, \medskip \\
U_1=\cdots U_{n-1}=0, \  \frac{\partial U_n}{\partial x_n}=k(x)  U_n \quad  &\mbox{on} \ \partial\mathbb B^+_1 \cap \{x|x_n=0\},
\end{array}
\right.
\label{system}
\end{equation}
where $A(x)$, $B(x)$ are matrix functions and $k(x)=- \frac{\partial  \ln g^{\frac{1}{2}}}{\partial x_n}$ is a scalar function which all depend on the metric $g$ and its derivatives.

 The rest of section is devoted to the statement of the quantitative Carleman estimates for elliptic systems in half balls with proper boundary conditions. The quantitative Carleman estimates are important tools for the study of three-ball inequalities and doubling inequalities in the next section.
Let $r=r(y)$ be the Riemannian distance from the origin to $y$, which is always less than the
injectivity radius.  Carleman estimates are weighted integral inequalities with some
 weight function $\exp\{\beta \psi(x)\}$. We construct the weight function $\psi$ as follows. Let $\psi(y)=-\phi(\ln r(x))$, where
 $\phi(t)=t+\ln t^2$ for $(-\infty, \ T_0]$ and $T_0$ is negative with $|T_0|$ is sufficiently large enough. It is easy to see that
  the function $\phi(t)$ satisfies the
following properties
\begin{align}
1+\frac{2}{T_0}\leq \phi'(t)\leq 1, \label{ass1}\\
\lim\limits_{t\to-\infty}\frac{-\phi''(t)}{e^t}=+\infty.
\label{ass2}
\end{align}
  We are able to establish  the following Carleman estimates. Note that $E(x)$ and $h(x)$ below may have large $C^1$ norm. The notation $\|\cdot\|_{\mathbb B^+_R}$ denotes the $L^2$ norm in the half ball $\mathbb B^+_R$ unless otherwise stated.

\begin{proposition} There exist positive constants $C_1$, $C_0$ and small constant $r_0$ such that for
any $V\in C^{\infty}_{0}\in \mathbb B^+_{r_0}\backslash \{0\}$, $E(x), h(x)\in C^1$, and $$\beta>C_0(1+\sqrt{\|E\|_{C^1}}+ \|h\|_{C^1}),$$ one has
\begin{align}
C_1\|r^2 e^{\beta \psi}F\|_{\mathbb B^+_{r_0}} \geq  \beta^\frac{3}{2} \|  e^{\beta \psi} (\log r)^{-1} V\|_{\mathbb B^+_{r_0}}
+\beta^\frac{1}{2} \| r e^{\beta \psi} (\log r)^{-1} \nabla V\|_{\mathbb B^+_{r_0}},
 \label{Carle}
\end{align}
 where the vector function $ V(x)$ satisfies
 \begin{align}
 \triangle V+E(x)V=F(x) \quad \quad \quad &\mbox{in} \ \mathbb B^+_{r_0}, \nonumber \\
 V_1=0, \cdots, V_{n-1}=0, \ \frac{\partial V_n}{\partial x_n}=h(x) V_n, \quad \quad &\mbox{on}\ \partial\mathbb B^+_{r_0}\cap\{x|x_n=0\}.
 \label{modelll}
 \end{align}
 Furthermore, if $V\in C^{\infty}_{0}\in \mathbb B^+_{r_0}\backslash \mathbb B^+_{\rho}$, then
 \begin{align}
C_1\|r^2 e^{\beta \psi}F\|_{\mathbb B^+_{r_0}} &\geq  \beta^\frac{3}{2} \|  e^{\beta \psi} (\log r)^{-1} V\|_{\mathbb B^+_{r_0}}+\beta \rho^\frac{1}{2} \| r^{-\frac{1}{2}} e^{\beta \psi}  V\|_{\mathbb B^+_{r_0}} \nonumber \\
&+\beta^\frac{1}{2} \| r e^{\beta \psi} (\log r)^{-1} \nabla V\|_{\mathbb B^+_{r_0}}.
 \label{Strongcarle}
\end{align}
 \label{pro2}
\end{proposition}

Next we show a lemma which will be applied for the elliptic system (\ref{system}) in our latter arguments.
\begin{lemma} There exist positive constants $C_1$, $C_0$ and small constant $r_0$ such that for
any $V\in C^{\infty}_{0}(\mathbb B^+_{r_0})\backslash \{0\}$, and $$\beta>C_0(1+\sqrt{\lambda}),$$ one has
\begin{align}
C_1\|r^2 e^{\beta \psi}F\|_{\mathbb B^+_{r_0}} \geq  \beta^\frac{3}{2} \|  e^{\beta \psi} (\log r)^{-1} V\|_{\mathbb B^+_{r_0}}
+\beta^\frac{1}{2} \| r e^{\beta \psi} (\log r)^{-1} \nabla V\|_{\mathbb B^+_{r_0}},
 \label{Carle2}
\end{align}
 where $ V(x)$ satisfies
 \begin{align}
 \triangle V+\langle  B, \nabla V\rangle+AV+\lambda V=F \quad \quad \quad &\mbox{in} \ \mathbb B^+_{r_0}, \nonumber \\
 V_1=0, \cdots, V_{n-1}=0, \ \frac{\partial V_n}{\partial x_n}=k(x) V_n, \quad \quad &\mbox{on}\ \partial\mathbb B^+_{r_0}\cap\{x|x_n=0\}.
 \label{model}
 \end{align}
 Furthermore, if $V\in C^{\infty}_{0}\in \mathbb B^+_{r_0}\backslash \mathbb B^+_{\rho}$, then
 \begin{align}
C_1\|r^2 e^{\beta \psi}F\|_{\mathbb B^+_{r_0}} &\geq  \beta^\frac{3}{2} \|  e^{\beta \psi} (\log r)^{-1} V\|_{\mathbb B^+_{r_0}}+\beta \rho^\frac{1}{2} \| r^{-\frac{1}{2}} e^{\beta \psi}  V\|_{\mathbb B^+_{r_0}} \nonumber \\
&+\beta^\frac{1}{2} \| r e^{\beta \psi} (\log r)^{-1} \nabla V\|_{\mathbb B^+_{r_0}}.
 \label{Strongcarle2}
\end{align}
\end{lemma}
\begin{proof}
We assume that $\beta>C_0(1+\sqrt{\lambda})$ and $\lambda$ is large enough. Since  $B(x), A(x)$ and $k(x)$ are fixed functions depending only on the metric $g$ and its derivatives, from triangle inequality and (\ref{Carle}), we readily have
\begin{align}
C_1\|r^2 e^{\beta \psi}(\triangle V+\langle  B, \nabla V\rangle+AV+\lambda V)\|_{\mathbb B^+_{r_0}} &\geq  \|r^2 e^{\beta \psi}(\triangle V+\lambda V)\|_{\mathbb B^+_{r_0}}-\|r^2 e^{\beta \psi} \langle  B, \nabla V\rangle\|_{\mathbb B^+_{r_0}} \nonumber \\
&-\|r^2 e^{\beta \psi}|AV|\|_{\mathbb B^+_{r_0}}\nonumber \\
&
\geq \beta^\frac{3}{2} \|  e^{\beta \psi} (\log r)^{-1} V\|_{\mathbb B^+_{r_0}} \nonumber \\&+\beta^\frac{1}{2} \| r e^{\beta \psi} (\log r)^{-1} \nabla V\|_{\mathbb B^+_{r_0}}
-\|r^2 e^{\beta \psi} \langle  B, \nabla V\rangle\|_{\mathbb B^+_{r_0}} \nonumber \\
&-\|r^2 e^{\beta \psi}|AV|\|_{\mathbb B^+_{r_0}}.
\label{lala}
\end{align}
By the assumption of $\beta$ and the fact that $r$ is small, the last two terms in the last inequality can be controlled by the right hand side of (\ref{Carle})). Thus, the Carleman estimates (\ref{Carle2}) can be derived from (\ref{lala}). The conclusion (\ref{Strongcarle2}) can be arrived using the same strategy from the estimates (\ref{Strongcarle}).
\end{proof}
Inspired by the approach in \cite{DF}, \cite{B}, \cite{Zh} and \cite{R}, we prove the quantitative Carleman estimates (\ref{Carle}) and (\ref{Strongcarle}) by introducing new coordinates and conjugate operators. Since the arguments of the Carleman estimates is long and detailed, we present the proof in the last section, i.e. Section 6. We will only use the Carleman estimates (\ref{Carle2}) and (\ref{Strongcarle2})  for the later proof of the doubling inequalities and the measure of critical sets.

\section{Three-ball inequalities and doubling inequalities}
In this section, thanks to the quantitative Carleman estimates (\ref{Carle2}) and (\ref{Strongcarle2}), we will first show the three half-ball inequalities of the gradient of Dirichlet eigenfunctions in the half balls centered at the boundary. Then we will obtain the quantitative doubling inequalities of the gradients in those half balls. In the proof of the doubling inequality, an important ingredient is a lower bound estimate for $L^2$ norm of the gradients of Dirichlet eigenfunctions, since the doubling manifold technique in \cite{DF} is not available directly, our new idea is to study the propagation between the neighborhood of the boundary and the interior of the manifold.

 Relied on the Carleman estimates (\ref{Carle2}), it is standard to
 have the following quantitative interior three-ball inequality
\begin{align}
\|\nabla e_\lambda \|_{L^\infty(\mathbb B_{\frac{R}{2}}(x))}\leq e^{C\sqrt{\lambda}} \|\nabla e_\lambda\|^{\hat{\tau}}_{L^\infty(\mathbb B_{\frac{R}{4}}(x))} \|\nabla e_\lambda\|^{1-\hat{\tau}}_{L^\infty(\mathbb B_{R}(x))}
\label{interiorball}
\end{align}
 if $\mathbb B_{R}(x)\subset \mathcal{M}$, where $0<\hat{\tau}<1$ is a constant.
 For the complete of the presentation, we present the approach to obtain these inequalities in half balls by considering the boundary conditions. Recall that $r(x)$ be the geodesic distance from $x$ to origin. By rotation and translation, we may assume that $0\in \partial \mathcal{M}$. Denote
$A_{R_1, R_2}=\{x\in \mathcal{M}| R_1\leq r(x)\leq R_2\}$ be the annulus. Let $\|U\|_{R_1, R_2}$ be the $L^2$ norm of $U$ on $A_{R_1, R_2}$.

\begin{lemma}
There exist positive constants $\bar R<r_0$, $C$ and $0<\tau<1$ which depend only on $\mathcal{M}$ such that, for any $R<\bar R$ and any $x_0\in \partial \mathcal{M}$, the Dirichlet eigenfunctions $e_\lambda$ of (\ref{diri}) satisfy
\begin{align}
\|\nabla e_\lambda \|_{L^\infty(\mathbb B^+_{2R}(x_0))}\leq e^{C\sqrt{\lambda}} \|\nabla e_\lambda\|^{\tau}_{L^\infty(\mathbb B^+_R(x_0))} \|\nabla e_\lambda\|^{1-\tau}_{L^\infty(\mathbb B^+_{3R}(x_0))}.
\label{three-balls}
\end{align}
\end{lemma}
\begin{proof}
We study the elliptic systems (\ref{system}) derived from the derivatives of Dirichlet eigenfunctions. By translation, let $x_0=0$.
To apply the Carleman estimates (\ref{Carle2}),
we introduce a smooth cut-off function $\phi(r)\in C^\infty_0(\mathbb B_{3R})$ with $R<\frac{r_0}{8}$. Let $0<\phi(r)<1$ satisfy the following properties:
\begin{itemize}
\item $\phi(r)=0$ \ \ \mbox{if} \ $r(x)<\frac{R}{5}$ \ \mbox{or} \  $r(x)>\frac{7R}{3}$, \medskip
\item $\phi(r)=1$ \ \ \mbox{if} \ $\frac{3R}{5}<r(x)<\frac{13R}{6}$, \medskip
\item $|\nabla^\alpha \phi|\leq \frac{C}{R^{|\alpha|}}$
\end{itemize}
for $\alpha=(\alpha_1, \cdots, \alpha_n)$. Thus, the function $\phi U$ is supported in the annulus $A_{\frac{R}{5}, \frac{7R}{3}} $. Applying the Carelman estimates (\ref{Carle2}) with $V=\phi U$ and considering the elliptic systems (\ref{system}), we obtain that
\begin{align}
\beta^\frac{3}{2} \|(\log r)^{-1} e^{\beta \psi}\phi U\|_{\mathbb B^+_{r_0}}& \leq C\| r^2 e^{\beta \psi}\big(\triangle (\phi U)+\langle B, \  \nabla (\phi U)\rangle +A\phi U+\lambda \phi U\big)\|_{\mathbb B^+_{r_0}} \nonumber  \\
&=C \| r^2 e^{\beta \psi}\big(\triangle \phi U+ 2\nabla \phi \nabla U  +\langle B, \ \nabla \phi U\rangle \big)\|_{\mathbb B^+_{r_0}}.
\end{align}
It follows from the properties of $\phi$ that
\begin{align*}
\| e^{\beta \psi}U \|_{\frac{3R}{5}, \frac{13R}{6}}&\leq C (\| e^{\beta \psi} U\|_{\frac{R}{5}, \frac{3R}{5}}+\| e^{\beta  \psi} U\|_{\frac{13R}{6}, \frac{7R}{3}} ) \\
&+ C( \| r e^{\beta \psi} \nabla U \|_{\frac{R}{5}, \frac{3R}{5}}+\| r e^{\beta \psi} \nabla U\|_{\frac{13R}{6}, \frac{7R}{3}}).
\end{align*}
Note that the weight function $\psi$ is radial and decreasing. We obtain that
\begin{align}
 \|e^{\beta \psi}  U\|_{\frac{3R}{5}, \frac{13R}{6}}&\leq C (e^{\beta \psi(\frac{R}{5})} \|   U\|_{\frac{R}{5}, \frac{3R}{5}}+e^{\beta \psi(\frac{13R}{6})} \|  U\|_{\frac{13R}{6}, \frac{7R}{3}} ) \nonumber \\
&+ C(e^{\beta \psi (\frac{R}{5})}\| r \nabla U\|_{\frac{R}{5}, \frac{3R}{5}}+e^{\beta \psi(\frac{13R}{6})}\| r \nabla  U\|_{\frac{13R}{6}, \frac{7R}{3}}).
\label{recall}
\end{align}

The following Caccioppoli inequality for the elliptic systems (\ref{system})
\begin{align}
\| \nabla U\|_{\mathbb B^+_{c_2R}} \leq \frac{C (\sqrt{\lambda}+1)}{R}\|U\|_{\mathbb B^+_{ c_1R}}
\label{hihcac}
\end{align}
 holds for all positive constants $0<c_2<c_1<1$, which can be proved by multiplying the elliptic systems (\ref{system}) by $\hat{\phi}^2U$ for some cut-off function $\hat{\phi}$ and using the trace inequality (\ref{trace}).
 Thus, the estimate (\ref{recall}) implies that
\begin{align}
\|U\|_{\frac{3R}{5}, 2R}\leq C\sqrt{\lambda} \big( e^{\beta(\psi(\frac{R}{5})-\psi(2R))} \|U\|_{\mathbb B^+_{R}}+  e^{\beta(\psi(\frac{13R}{6})-\psi(2R))} \|U\|_{\mathbb B^+_{3R}}\big).
\label{inequal}
\end{align}
We choose some new parameters
$$ \tau^1_R=\psi(\frac{R}{5})-\psi(2R),$$
$$ \tau^2_R=\psi(2R)-\psi(\frac{13R}{6}).$$
Thanks to the definition of $\psi$, we learn that
$$ 0<\tau^{-1}_1<\tau^1_R<\tau_1 \quad \mbox{and} \quad 0<\tau_2<\tau^2_R<\tau^{-1}_2,$$
where $\tau_1$ and $\tau_2$ are independent of $R$.
Adding $\|U\|_{\frac{3R}{5}}$ to both sides of the inequality (\ref{inequal}) yields that
\begin{align}
\|U\|_{\mathbb B^+_{2R}}\leq C\sqrt{\lambda}\big( e^{\beta\tau_1}\|U\|_{\mathbb B^+_{R}}+  e^{-\beta\tau_2}\|u\|_{\mathbb B^+_{3R}}     \big).
\end{align}
We want to incorporate the second term in the right hand side of the last inequality into the left hand side. To this end,
 we choose $\beta$ such that
$$C\sqrt{\lambda} e^{-\beta\tau_2}\|U\|_{\mathbb B^+_{3R}}\leq \frac{1}{2}\|U\|_{\mathbb B^+_{2R}},   $$
which holds if
$$\beta\geq \frac{1}{\tau_2} \ln \frac{2C\sqrt{\lambda}\|U\|_{\mathbb B^+_{3R}}}{\|U\|_{\mathbb B^+_{2R}} }.   $$
Thus, we derive that
\begin{equation}
\|U\|_{\mathbb B^+_{2R}}\leq C\sqrt{\lambda} e^{\beta\tau_1}\|U\|_{\mathbb B^+_{R}}.
\label{substitute}
\end{equation}
Recall that the lower bound $\beta>C\sqrt{\lambda}$ is required to apply the Carleman estimates (\ref{Carle2}). Hence we choose
$$ \beta=C \sqrt{\lambda}+\frac{1}{\tau_2} \ln \frac{2C\sqrt{\lambda}\|U\|_{\mathbb B^+_{3R}}}{\|U\|_{\mathbb B^+_{2R}} }. $$
Substituting such $\beta$ in (\ref{substitute}) gives that
\begin{align}
\|U\|_{\mathbb B^+_{2R}}^{\frac{\tau_2+\tau_1}{\tau_2}} \leq e^{ C\sqrt{\lambda}}\|U\|_{\mathbb B^+_{3R}}^{\frac{\tau_1}{\tau_2}} \|U\|_{\mathbb B^+_{R}}.
\end{align}
Raising the exponent $\frac{\tau_2}{\tau_2+\tau_1}$ to both sides of the last inequality yields that
\begin{align}
\|U\|_{\mathbb B^+_{2R}} \leq e^{ C\sqrt{\lambda}}\|U\|_{\mathbb B^+_{3R}}^{\frac{\tau_1}{\tau_1+\tau_2}} \|U\|_{\mathbb B^+_{R}}^{\frac{\tau_2}{\tau_1+\tau_2}}.
\end{align}
Set $\tau={\frac{\tau_2}{\tau_1+\tau_2}}$. Thus,  $0<\tau<1$.  Recall that $U=(\frac{\partial e_\lambda}{\partial x_1}, \cdots, \frac{\partial e_\lambda}{\partial x_n})$, we arrive at
\begin{align}
\|\nabla e_\lambda\|_{\mathbb B^+_{2R}} \leq e^{ C\sqrt{\lambda}} \|\nabla e_\lambda\|_{\mathbb B^+_{R}}^{\tau} \|\nabla e_\lambda\|_{\mathbb B^+_{3R}}^{1-\tau }.
\label{aimm}
\end{align}
By the standard elliptic estimates, the $L^\infty$ norm and $L^2$ norm of $U$ are comparable in different sizes of balls. We save the efforts in changing the sizes of balls in converting the $L^2$ norm in (\ref{aimm}) into the $L^\infty$ norm. 
 Thus,
the three half-ball inequality in the lemma is completed.
\end{proof}

Next we aim to find some quantitative lower bounds for $L^2$ norm of $\nabla e_\lambda$ in half-ball centered at the boundary $\partial\mathcal{M}$. Unlike the manifold without boundary or double manifold, we can not iterate the three-ball inequality directly to achieve the goal. Instead, we combine the three half-ball inequality  and interior three-ball inequality to do the propagation between the neighborhood  on the boundary and the interior of the manifold.

First we claim that
\begin{align}
\|\nabla e_\lambda\|_{L^\infty(\mathcal{M}_{R})} \geq e^{-C(R)\sqrt{\lambda}} \|\nabla e_\lambda\|_{L^\infty(\mathcal{M})},
\label{claim}
\end{align}
where $C(R)$ is a positive constant depending on $R$ and $\mathcal{M}_{R}=\{x\in \mathcal{M}| dist(x, \ \partial\mathcal{M})\leq {R}\}.$
Since $\nabla e_\lambda$ is smooth, there exist some point $\hat{x}\in \overline{\mathcal{M}}$ such that $|\nabla e_\lambda(\hat{x})|=\|\nabla e_\lambda\|_{L^\infty(\mathcal{M})}$. We may rescale to assume that $\|\nabla e_\lambda\|_{L^\infty(\mathcal{M})}=1$.
If $\hat{x}\in \mathcal{M}_{R}$, the claim (\ref{claim}) follows immediately. If $\hat{x}\in \mathcal{M}\backslash\mathcal{M}_{R}$,
we apply the interior three-ball inequality and three half-ball inequality. Assume that there exists some point $\hat{x}_0\in \mathcal{M}_{R}$ such that $|\nabla e_\lambda(\hat{x}_0)|=\|\nabla e_\lambda\|_{L^\infty(\mathcal{M}_{R})}$. Then we choose some point $\bar x\in \partial\mathcal{M}$ such that $ \hat{x}_0\in \mathbb B^+_R(\bar x)$ and
\begin{align}\|\nabla e_\lambda\|_{L^\infty(\mathbb B^+_{R}(\bar x))}=\|\nabla e_\lambda\|_{L^\infty(\mathcal{M}_{R})}=\delta_0
\label{dddelta}
\end{align}
for some $0<\delta_0<1$.
Applying the three half-ball inequality (\ref{three-balls}) at $\bar x$, we get
\begin{align}
\|\nabla e_\lambda \|_{L^\infty(\mathbb B^+_{2R}(\bar x))}\leq e^{C\sqrt{\lambda}} \delta_0^{\tau}.
\end{align}
Then we choose a point $x_1 \in \mathbb B^+_{2R}(\bar x)$ such that
\begin{align}
\mathbb B_{\frac{R}{4}}( x_1) \subset \mathbb B^+_{2R}(\bar x)\quad \mbox{and}  \quad \mathbb B_{\frac{R}{2}}( x_1) \not \subset\mathbb B^+_{2R}(\bar x).
\end{align}
Applying the interior three-ball inequality at $x_1$, we have
\begin{align}
\|\nabla e_\lambda \|_{L^\infty(\mathbb B_{\frac{R}{2}}(x_1))}&\leq e^{C\sqrt{\lambda}} \|\nabla e_\lambda\|^{\hat{\tau}}_{L^\infty(\mathbb B_{\frac{R}{4}}(x_1))} \|\nabla e_\lambda\|^{1-\hat{\tau}}_{L^\infty(\mathbb B_{R}(x_1))} \nonumber \\
&\leq e^{C\sqrt{\lambda}(1+\hat{\tau})} \delta_0^{\tau \hat{\tau}}.
\end{align}
Fix such  $R$, we choose a sequence of balls $\mathbb B_{\frac{R}{4}}(x_i)$ centered at $x_i$ such that $x_{i+1}\in \mathbb  B_{\frac{R}{4}}(x_i)$ and $B_{\frac{R}{4}}(x_{i+1})\subset \mathbb B_{\frac{R}{2}}(x_i)$. After finitely many of steps, we could get to the point $\hat{ x}$ where $|\nabla e_\lambda|(\hat {x})=1$, that is, $ x_1, x_2, \cdots, x_m=\hat {x}$. The number of $m$ depends on $R$ and $diam(\mathcal{M})$. Repeating the three-ball inequality (\ref{interiorball}) at those $x_i$, $i=2, 3, \cdots, m$, we arrive at
\begin{align}
\|\nabla e_\lambda\|_{L^\infty(\mathbb B_{\frac{R}{4}}(x_m))}\leq e^{C\sqrt{\lambda}(1+\hat{\tau}+\hat{\tau}^2+\cdots+\hat{\tau}^m)   } \delta_0^{\tau \hat{\tau}^m}.
\end{align}
Since $0<\tau, \hat{\tau}<1$, we obtain that
\begin{align}
\|\nabla e_\lambda\|_{L^\infty(\mathcal{M}_{R})}&\geq e^{\frac{-C\sqrt{\lambda}(1+\hat{\tau}+\hat{\tau}^2+\cdots+\hat{\tau}^m)}{\tau \hat{\tau}^m }}\nonumber \\
&\geq e^{-{C(R){\sqrt{\lambda}}}}\|\nabla e_\lambda\|_{L^\infty(\mathcal{M})}.
\end{align}
This verifies the claim (\ref{claim}).

Next we do a propogation of smallness using the three half-ball inequality to get a $L^\infty$ lower bound of $\nabla e_\lambda$ in the half-ball. Let $\hat{x}_1$ be any point on $\partial\mathcal{M}$. The application of three half-ball inequality (\ref{three-balls}) yields that
\begin{align}
\|\nabla e_\lambda \|_{L^\infty(\mathbb B^+_{2R}(\hat{x}_1))}&\leq e^{C\sqrt{\lambda}} \|\nabla e_\lambda\|^{{\tau}}_{L^\infty(\mathbb B^+_{R}(\hat{x}_1))}
\end{align}
since we have assumed that $\|\nabla e_\lambda \|_{L^\infty(\mathcal{M})}=1$. We choose $\hat{x}_2\in \partial\mathcal{M}$ such that
$\mathbb B^+_{R}(\hat{x}_2)\subset \mathbb B^+_{2R}(\hat{x}_1)$. Thus, we have
\begin{align}
\|\nabla e_\lambda \|_{L^\infty(\mathbb B^+_{R}(\hat{x}_2))}&\leq e^{C\sqrt{\lambda}} \|\nabla e_\lambda\|^{{\tau}}_{L^\infty(\mathbb B^+_{R}(\hat{x}_1))}
\end{align}
For such fixed $R$, we again choose a sequence of balls $\mathbb B^+_{R}(\hat{x}_i)$ centered at $\hat{x}_i\in \partial \mathcal{M}$ such that $\mathbb B^+_{R}(\hat{x}_{i+1})\subset \mathbb B^+_{2R}(\hat{x}_i)$. Since $\partial \mathcal{M}$ is compact, after finitely many of steps, we could get to the point $\bar x$. Recall the assumption of $\bar x$ in (\ref{dddelta}). That is, we choose a sequence of points, $\hat{x}_{1}, \hat{x}_{2}, \cdots, \hat{x}_{\bar m}=\bar x$. The number of $\bar m$ depends on $R$ and $\partial \mathcal{M}$.
 Repeating the three half-ball inequality (\ref{three-balls}) at those $\hat{x}_i$, $i=2, 3, \cdots, \bar m$, we arrive at
 \begin{align}
\|\nabla e_\lambda \|_{L^\infty(\mathbb B^+_{R}(\bar {x}))}&\leq e^{C\sqrt{\lambda}(1+\tau^1+\cdots+\tau^{\bar m-1})} \|\nabla e_\lambda\|^{{\tau^{\bar m}}}_{L^\infty(\mathbb B^+_{R}(\hat{x}_1))}.
\end{align}
Taking (\ref{claim}) and the assumption of $\bar x$ in (\ref{dddelta}) into consideration gives that
\begin{align}
\|\nabla e_\lambda \|_{L^\infty(\mathbb B^+_{R}(\hat{x}_1))}&\geq e^{\frac{-C\sqrt{\lambda}(1+\tau^1+\cdots+\tau^{\bar m-1})}{\tau^{\bar m}}} \|\nabla e_\lambda\|_{L^\infty(\mathcal{M}_R)} \nonumber \\
&\geq e^{-C_2(R)\sqrt{\lambda} } \|\nabla e_\lambda\|_{L^\infty(\mathcal{M})}.
\label{quant}
\end{align}

By rescaling, it also holds that
\begin{align}
\|\nabla e_\lambda \|_{L^\infty(\mathbb B^+_{\frac{R}{4}}(\hat{x}_1))}
\geq e^{-C_3(R)\sqrt{\lambda} } \|\nabla e_\lambda\|_{L^\infty(\mathcal{M})}.
\label{quantna}
\end{align}

Recall that the annulus $ A_{R_1, \ R_2}({x_0})=\{x\in \mathcal{M} | R_1\leq |x-x_0|\leq R_2\}$.
For any $x_0\in \partial\mathcal{M}$, there exist some point $\hat{x}_1\in \partial\mathcal{M}$ such that $\mathbb B^+_{\frac{R}{4}}(\hat{x}_1) \subset A_{\frac{R}{2}, \ R}(x_0)$. Therefore, (\ref{quantna}) also implies that
\begin{align}
\|\nabla e_\lambda\|_{L^\infty( A_{\frac{R}{2}, \ R}({x_0}))}
&\geq e^{-{C(R){\sqrt{\lambda}}}} \|\nabla e_\lambda\|_{L^\infty(\mathcal{M} )}.
\label{exdouex}
\end{align}

With aid of the quantitative lower bound (\ref{exdouex}) and the Carleman estimates  (\ref{Strongcarle2}),
we proceed to show the doubling inequality for the gradients of Dirichlet eigenfunctions $\nabla e_\lambda$ in the half ball.

\begin{proof}[Proof of Theorem \ref{th3}] We study the elliptic systems (\ref{system}) for $U$. Without loss of generality, let $x_0=0$.
Let us fix $R=\frac{\bar R}{8}$, where $\bar R$ is the one in the three half-ball inequality (\ref{three-balls}). Choose $0<\rho<\frac{R}{24}$ to be arbitrarily small. We introduce a smooth cut-off function $0<\phi<1$ as follows,
\begin{itemize}
\item $\phi(r)=0$ \ \ \mbox{if} \ $r(x)<\rho$ \ \mbox{or} \  $r(x)>2R$, \medskip
\item $\phi(r)=1$ \ \ \mbox{if} \ $2\rho<r(x)<\frac{3R}{2}$, \medskip
\item $|\nabla^\alpha \phi|\leq \frac{C}{\rho^\alpha}$ \ \ \mbox{if} $\rho<r(x)<2\rho$,\medskip
\item $|\nabla^\alpha \phi|\leq C$ \ \  \mbox{if} \ $\frac{3R}{2}<r(x)<2R$.
\end{itemize}
We apply the stronger Carleman estimates (\ref{Strongcarle2}) this time. Replacing $V$ by $\phi U$ and substituting it into (\ref{Strongcarle2}) gives that
\begin{align*}
 \| (\log r)^{-1} e^{\beta \psi} \phi U\|_{\mathbb B^+_{r_0}}&+ \beta\rho^\frac{1}{2} \| r^{\frac{-1}{2}} e^{\beta \psi}\phi U\|_{\mathbb B^+_{r_0}} \nonumber\\
&\leq C\| r^2 e^{\beta \psi}\big( \triangle (\phi U)+ \langle B, \ \ \nabla (\phi U)\rangle +A \cdot \nabla (\phi U)+\lambda \phi U\|_{\mathbb B^+_{r_0}}.
\end{align*}
It follows from the properties of $\phi$ and the elliptic systems (\ref{system})  that
\begin{align*}
 \| (\log r)^{-1} e^{\beta \psi} U\|_{\frac{R}{2}, R}+  \| e^{\beta \psi} U\|_{2\rho, 6\rho}
 &\leq C (\| e^{\beta \psi} U\|_{\rho, 2\rho}+\| e^{\beta \psi} U\|_{\frac{3R}{2}, 2R} ) \\
&+ C\| r e^{\beta \psi} \nabla U\|_{\rho, 2\rho}+\| r e^{\beta \psi} \nabla  U\|_{\frac{3R}{2}, 2R}).
\end{align*}
Note that $R$ is fixed. We take the exponential function $e^{\beta \psi}$ out of the norms by using the fact that $\psi$ is radial and decreasing. Thus,  we arrive at
\begin{align*}
e^{\beta \psi(R)}\|  U\|_{\frac{R}{2}, R}+ e^{\beta \psi({6\rho})}  \|  U\|_{2\rho, 6\rho}
&\leq C (e^{\beta \psi(\rho) }\|  U\|_{\rho, 2\rho}+e^{\beta\psi(\frac{3R}{2}) }\|  U\|_{\frac{3R}{2}, 2R} ) \\
&+ C(e^{\beta \psi(\rho)} \| r\nabla U\|_{\rho,
2\rho}+e^{\beta\psi(\frac{3R}{2})}\|r  \nabla U\|_{\frac{3R}{2}, 2R}).
\end{align*}
The application of  Caccioppoli  inequality (\ref{hihcac})
further implies that
\begin{align}
e^{\beta \psi(R)}\|  u\|_{\frac{R}{2}, R}+ e^{\beta \psi({6\rho})}  \|  u\|_{2\rho, 6\rho}
&\leq C \sqrt{\lambda} (e^{\beta \psi(\rho) }\|  U\|_{\mathbb B^+_{3\rho}}+e^{\beta\psi(\frac{3R}{2})}\|  U\|_{\mathbb B^+_{3R}}).
\end{align}
Adding $ e^{\beta \psi({6\rho})} \|U\|_{2\rho}$ to both sides of last inequality, we get that
\begin{align}
e^{\beta \psi(R)}\|  U\|_{\frac{R}{2}, R}+ e^{\beta \psi({6\rho})}  \|  U\|_{\mathbb B^+_{6\rho}}
&\leq C \sqrt{\lambda} (e^{\beta \psi(\rho) }\|  U\|_{\mathbb B^+_{3\rho}}+e^{\beta \psi(\frac{3R}{2})}\|  U\|_{\mathbb B^+_{3R}}).
\label{drop}
\end{align}
We want to remove the second term in the right hand side of the last inequality. To this end,
we choose $\beta$ to satisfy
$$C\sqrt{\lambda} e^{\beta\psi(\frac{3R}{2})}\|U\|_{\mathbb B^+_{3R}}\leq \frac{1}{2}e^{\beta\psi(R)} \|U\|_{\frac{R}{2}, R}.  $$
That is, at least
$$ \beta\geq \frac{1}{\psi(R)-\psi(\frac{3R}{2})}\ln \frac{ 2C\sqrt{\lambda} \|U\|_{\mathbb B^+_{3R}}}{ \|U\|_{\frac{R}{2}, R}}   $$
is needed.
Thus, we obtain that
\begin{align}
e^{\beta \psi(R)}\|  U\|_{\frac{R}{2}, R}+ e^{\beta \psi({6\rho})}  \|  U\|_{\mathbb B^+_{6\rho}}
\leq C \sqrt{\lambda} e^{\beta \psi(\rho) }\|  U\|_{\mathbb B^+_{3\rho}}.
\label{droppd}
\end{align}

To apply the Carleman estimates (\ref{Strongcarle2}), we need the assumption that  $\beta\geq C\sqrt{\lambda}$. Therefore, we select
$$\beta=C\sqrt{\lambda}+ \frac{1}{\psi(R)-\psi(\frac{3R}{2})}\ln \frac{ 2C\sqrt{\lambda} \|U\|_{\mathbb B^+_{3R}}}{ \|U\|_{\frac{R}{2}, R}}. $$
Furthermore, dropping the first term in (\ref{droppd}), we derive that
\begin{align}
\|U\|_{\mathbb B^+_{6\rho}}&\leq C\sqrt{\lambda} \exp\{ \big(C\sqrt{\lambda}+ \frac{1}{\psi(R)-\psi(\frac{3R}{2})}\ln \frac{ 2C\sqrt{\lambda} \|U\|_{\mathbb B^+_{3R}}}{ \|U\|_{\frac{R}{2}, R}}\big)\big(\psi(\rho)-\psi(6\rho)\big)  \}\|U\|_{\mathbb B^+_{3\rho}} \nonumber \\
&\leq e^{C\sqrt{\lambda}}   (\frac{\|U\|_{\mathbb B^+_{3R}}}{ \|U\|_{\frac{R}{2}, R}})^C \|U\|_{\mathbb B^+_{3\rho}},
\label{tata}
\end{align}
where we have used the fact that the bounds  $\psi(R)-\psi(\frac{3R}{2})$ and $\psi(\rho)-\psi(6\rho)$ are independent of $R$ or $\rho$.
Since $U=(\frac{\partial e_\lambda}{\partial x_1}, \cdots, \frac{\partial e_\lambda}{\partial x_n})$, it follows from (\ref{exdouex}) that
$$\frac{\|\nabla e_\lambda\|_{\mathbb B^+_{3R}}}{ \|\nabla e_\lambda\|_{\frac{R}{2}, R}}\leq e^{C\sqrt{\lambda}}. $$
Together  the last inequality with (\ref{tata}), we derive that
$$ \|\nabla e_\lambda\|_{\mathbb B^+_{6\rho}}\leq  e^{C\sqrt{\lambda}} \|\nabla e_\lambda\|_{\mathbb B^+_{3\rho}}. $$
Let $\rho=\frac{r}{3}$. The doubling inequality
\begin{align}
\|\nabla e_\lambda\|_{\mathbb B^+_{2r}}\leq  e^{C\sqrt{\lambda}} \|\nabla e_\lambda\|_{\mathbb B^+_{r}}
\end{align}
follows for $r\leq \frac{R}{8}$. If $r\geq \frac{R}{8}$, since $R$ is fixed, from (\ref{quant}), we can also derive that
\begin{align}
\|\nabla e_\lambda\|_{\mathbb B^+_{2r}}\leq  e^{C\sqrt{\lambda}} \|\nabla e_\lambda\|_{\mathbb B^+_{r}}.
\label{niuniu}
\end{align}
Thus, we have obtained (\ref{niuniu})
for any $0<r<r_0$, where $C$ only depends on the manifold $\mathcal{M}$. Since
 the $L^\infty$ norm and $L^2$ norm are equivalent for the elliptic equations, we arrive at the conclusions in  the theorem.
\end{proof}

To obtain the bounds of interior critical sets, we will make use of interior doubling inequalities. We do the propogation arguments between the consecutive annulus of the manifold.
\begin{lemma}
Let $e_\lambda$ be the Dirichlet eigenfunctions (\ref{diri}). For any fixed $R>0$, there exists a positive constant $C$ depending only on the manifold $\mathcal{M}$ and $R$ such that
\begin{align}
\|\nabla e_\lambda \|_{L^\infty(\mathbb B_{2r}(x_0))}\leq e^{ C\sqrt{\lambda} }\|\nabla e_\lambda\|_{L^\infty(\mathbb B_{r}(x_0))}
\label{III2}
\end{align}
for any $x_0\in \mathcal{M}\backslash \mathcal{M}_{2R}$ and $0<r\leq \frac{R}{2}$.
\label{submit}
\end{lemma}
\begin{proof}
We first claim that
\begin{align}
\|\nabla e_\lambda\|_{L^\infty(\mathcal{M}\backslash\mathcal{M}_{2R})} \geq e^{-C(R)\sqrt{\lambda}} \|\nabla e_\lambda\|_{L^\infty(\mathcal{M}\backslash\mathcal{M}_{R})}.
\label{claim2}
\end{align}
Let us normalize $\|\nabla e_\lambda\|_{L^\infty(\mathcal{M}\backslash\mathcal{M}_{R})}=1$. Assume that there exists some point $\hat{x}\in \mathcal{M}\backslash\mathcal{M}_{R}$ such that $\|\nabla e_\lambda\|_{L^\infty(\mathcal{M}\backslash\mathcal{M}_{R})}=|\nabla e_\lambda(\hat{x})|$. If $\hat{x}\in \mathcal{M}\backslash\mathcal{M}_{2R}$, the claim (\ref{claim2}) holds directly. If $\hat{x}\in \mathcal{M}_{2R}\backslash\mathcal{M}_{R}$, we will use the interior three-ball inequality to do the propagation. Let $\| \nabla e_\lambda\|_{L^\infty(\mathcal{M}\backslash\mathcal{M}_{2R})}=\bar \delta<1$. There exists some point $\hat{x}_0$ such that $\| \nabla e_\lambda\|_{L^\infty(\mathcal{M}\backslash\mathcal{M}_{2R})}=| \nabla e_\lambda( \hat{x}_0)|$.  Applying the three-ball inequality (\ref{interiorball}) at $\hat{x}_0$ gives that
\begin{align}
\|\nabla e_\lambda \|_{L^\infty(\mathbb B_{\frac{R}{2}}(\hat{x}_0))}&\leq e^{C\sqrt{\lambda}} \|\nabla e_\lambda\|^{\hat{\tau}}_{L^\infty(\mathbb B_{\frac{R}{4}}(\hat{x}_0))} \nonumber  \\
&\leq e^{C\sqrt{\lambda}} \bar \delta^{\hat{\tau}}.
\end{align}

Fix such $R$, we choose a sequence of balls $\mathbb B_{\frac{R}{4}}(x_i)$ centered at $x_i$ such that $x_{i+1}\in \mathbb  B_{\frac{R}{4}}(x_i)$ and $B_{\frac{R}{4}}(x_{i+1})\subset \mathbb B_{\frac{R}{2}}(x_i)$. After finitely many of steps, we could get to the point $\hat{ x}$ where $|\nabla e_\lambda|(\hat {x})=1$, that is, $ \hat{x}_0=x_1, x_2, \cdots, x_m=\hat {x}$. The number of $m$ depends on $R$ and $diam(\mathcal{M})$. Repeating the three-ball inequality (\ref{interiorball}) at those $x_i$, $i=1, 3, \cdots, m$, we arrive at
\begin{align}
\|\nabla e_\lambda\|_{L^\infty(\mathbb B_{\frac{R}{4}}(x_m))}\leq e^{C\sqrt{\lambda}(1+\hat{\tau}+\hat{\tau}^2+\cdots+\hat{\tau}^{m-1})   } \bar\delta^{ \hat{\tau}^m}.
\end{align}
Since $0<\hat{\tau}<1$, we obtain that
\begin{align}
\|\nabla e_\lambda\|_{L^\infty(\mathcal{M}\backslash\mathcal{M}_{2R})}&\geq e^{\frac{-C\sqrt{\lambda}(1+\hat{\tau}+\hat{\tau}^2+\cdots+\hat{\tau}^m)}{ \hat{\tau}^m }}\nonumber \\
&\geq e^{-{C(R){\sqrt{\lambda}}}}\|\nabla e_\lambda\|_{L^\infty(\mathcal{M}\backslash\mathcal{M}_{R})}.
\end{align}
This verifies the claim (\ref{claim2}).

Next we do a propogation of smallness using the three-ball inequality (\ref{interiorball}) to get a $L^\infty$ lower bound of $\nabla e_\lambda$ in $\mathcal{M}\backslash\mathcal{M}_{2R}$.  Choosing any point $\hat{x}_1\in \mathcal{M}\backslash\mathcal{M}_{2R}$, we apply the three-ball inequality (\ref{interiorball}) to have
\begin{align}
\|\nabla e_\lambda \|_{L^\infty(\mathbb B_{\frac{R}{2}}(\hat{x}_1))}&\leq e^{C\sqrt{\lambda}} \|\nabla e_\lambda\|^{{\hat{\tau}}}_{L^\infty(\mathbb B_{\frac{R}{4}}(\hat{x}_1))}
\end{align}
since we have assumed that $\|\nabla e_\lambda \|_{L^\infty(\mathcal{M}\backslash\mathcal{M}_{R})}=1$. We choose $\hat{x}_2\in \mathcal{M}\backslash\mathcal{M}_{2R}$ such that
$\mathbb B_{\frac{R}{4}}(\hat{x}_2)\subset \mathbb B_{\frac{R}{2}}(\hat{x}_1)$. Hence, it implies that
\begin{align}
\|\nabla e_\lambda \|_{L^\infty(\mathbb B_{\frac{R}{4}}(\hat{x}_2))}&\leq e^{C\sqrt{\lambda}} \|\nabla e_\lambda\|^{{\hat{\tau}}}_{L^\infty(\mathbb B_{\frac{R}{4}}(\hat{x}_1))}.
\end{align}
For such $R$, we again choose a sequence of balls $\mathbb B_{\frac{R}{4}}(\hat{x}_i)$ centered at $\hat{x}_i\in \mathcal{M}\backslash \mathcal{M}_{2R}$ such that $\mathbb B_{\frac{R}{4}}(\hat{x}_{i+1})\subset \mathbb B_{\frac{R}{2}}(\hat{x}_i)$. Since the closure of $ \mathcal{M}\backslash\mathcal{M}_{2R}$ is compact,  finitely many of steps of iterations lead to the point $\hat{x}_0$ where the maximum of $|\nabla e_\lambda|$ is achieved in $\mathcal{M}\backslash\mathcal{M}_{2R}$. That is, we choose a sequence of points, $\hat{x}_{1}, \hat{x}_{2}, \cdots, \hat{x}_{\bar m}=\hat{x}_0$. The number of $\bar m$ depends on $R$ and $\mathcal{M}$.
 Repeating the three-ball inequality (\ref{interiorball}) at those $\hat{x}_i$, $i=1, 2, 3, \cdots, \bar m$, we arrive at
 \begin{align}
\|\nabla e_\lambda \|_{L^\infty(\mathbb B_{\frac{R}{4}}(\hat{x}_0))}&\leq e^{C\sqrt{\lambda}(1+\hat{\tau}^1+\cdots+\hat{\tau}^{\bar m-1})} \|\nabla e_\lambda\|^{{\hat{\tau}^{\bar m}}}_{L^\infty(\mathbb B_{\frac{R}{4}}(\hat{x}_1))}.
\end{align}
Taking (\ref{claim2})  into consideration gives that
\begin{align}
\|\nabla e_\lambda \|_{L^\infty(\mathbb B_{\frac{R}{4}}(\hat{x}_1))}&\geq e^{\frac{-C\sqrt{\lambda}(1+\hat{\tau}^1+\cdots+\hat{\tau}^{\bar m-1})}{\hat{\tau}^{\bar m}}} \|\nabla e_\lambda\|_{L^\infty(\mathcal{M} \backslash\mathcal{M}_{2R} )} \nonumber \\
&\geq e^{-C_2(R)\sqrt{\lambda} } \|\nabla e_\lambda\|_{L^\infty(\mathcal{M}\backslash\mathcal{M}_R )}.
\label{quantfu}
\end{align}

Recall that the annulus $ A_{\frac{R}{4}, \ \frac{R}{2}}({x_0})=\{x\in \mathcal{M} | \frac{R}{4}\leq |x-x_0|\leq \frac{R}{2}\}$.
For any $x_0\in \mathcal{M} \backslash\mathcal{M}_{2R}$, there exist some point $\hat{x}_1\in \mathcal{M} \backslash\mathcal{M}_{2R}$ such that $\mathbb B_{\frac{R}{4}}(\hat{x}_1) \subset A_{\frac{R}{4}, \ \frac{3R}{4}}(x_0)$. Therefore, it follows from (\ref{quantfu}) that
\begin{align}
\|\nabla e_\lambda\|_{L^\infty( A_{\frac{R}{4}, \ \frac{3R}{4}}({x_0}))}
&\geq e^{-{C(R){\sqrt{\lambda}}}} \|\nabla e_\lambda\|_{L^\infty(\mathcal{M}\backslash\mathcal{M}_R  )}.
\label{exdouex22}
\end{align}

Following the proof of Theorem \ref{th3} and using the estimates (\ref{exdouex22}), we are able to obtain the interior doubling inequality (\ref{III2}).
\end{proof}

\section{Upper bounds of critical sets}

In this section, we will prove the upper bounds for the interior critical sets of Dirichlet eigenfunctions in real analytic manifolds. We first prove the upper bounds in the neighborhood of the boundary $\partial\mathcal{M}$. Then we show the upper bounds in the interior of the manifold $\mathcal{M}$. We will make use of doubling inequalities (\ref{LLL}), (\ref{III2}), and Lemma \ref{wwhy}. The main idea is similar to the proof of Theorem \ref{th1}. For the complete of the presentation, we provide the details.

 For the upper bounds in the neighborhood of the boundary, we  study the gradient of Dirichlet eigenfunctions in (\ref{system}) in an extended region. From the analyticity results in \cite{MN} or \cite{Mo}, it is clear that $U(x)=(\frac{\partial e_\lambda}{\partial x_1}, \frac{\partial e_\lambda}{\partial x_2},\cdots, \frac{\partial e_\lambda}{\partial x_n})$ is real analytic if the manifold $(\mathcal{M}, g)$ is analytic.
To do the analytic continuation across the boundary, we get rid of the large parameter $\lambda$. We introduce the following lifting argument.
Let \begin{align}\hat{U}(x, t)=e^{\sqrt{\lambda} t} U(x).\label{udefi}\end{align}

Since $U(x)$ satisfies the system of equations (\ref{system}), then $\hat{U}(x, t)$ satisfies the equation
\begin{equation}
\left \{ \begin{array}{rll}
\triangle_g \hat{U}+\partial_t^2 \hat{U} +\langle B, \  \nabla \hat{U}\rangle +A \cdot \hat{U}=0 \quad &\mbox{in} \ \mathcal{M}\times (-\infty, \ -\infty),  \medskip\\
\hat{U}_1=\cdots \hat{U}_{n-1}=0, \  \frac{\partial \hat{U}_n}{\partial x_n}=k(x)  \hat{U}_n  \quad &\mbox{on} \ {\partial\mathcal{M}}\times (-\infty, \ -\infty).
\end{array}
\right.
\label{killz}
\end{equation}

 We
introduce the ball with  as
\begin{align*}
\Omega_{R}=\{ (x, t)\in \mathbb R^{n+1}| |x|<R,  \ |t|<R\}
\end{align*}
and half-cube
\begin{align*}
\Omega^+_{R}=\{ (x, t)\in \mathbb R^{n+1} | |x|<R \ \mbox{with} \  x_n\geq 0, \ |t|<R\}.
\end{align*}
Choose any point $p\in \partial\mathcal{M}$, using Fermi geodesic coordinates and rescaling arguments, we may study the function $ \hat{U}(x, t )$  locally in the ball centered at origin with the flatten boundary. Hence, $ \hat{U}(x, t )$ satisfies the following equation locally
\begin{equation}
\left \{ \begin{array}{lll}
\triangle_g \hat{U}+\partial_t^2 \hat{U} +\langle B, \  \nabla \hat{U}\rangle +A \cdot \hat{U}=0 \quad &\mbox{in} \ \Omega^+_2,  \medskip\\
\hat{U}_1=\cdots \hat{U}_{n-1}=0, \  \frac{\partial \hat{U}_n}{\partial x_n}=k(x)  \hat{U}_n \quad &\mbox{on} \  \Omega^+_{2}\cap\{x_n=0\}.
\end{array}
\right.
\label{backgo}
\end{equation}

 By the analyticity results in \cite{MN} or \cite{Mo}, we can extend $\hat{U}(x, t)$ to the region $\Omega_{\rho}$, where $\rho>0$ depends only on $\mathcal{M}$. Moreover, we have the following growth control estimates
\begin{align}
\|\hat{U}\|_{L^\infty(\Omega_{\rho})}
\leq C\|\hat{U}\|_{L^\infty (\Omega^+_{2})},
\label{agree}
\end{align}
where $C$ depends only on $\mathcal{M}$. Since $\mathcal{M}$ is a real analytic Riemannian manifold with
boundary, we may embed $\mathcal{M}\subset \mathcal{M}_1$ as a
relatively compact subset, where $\mathcal{M}_1$ is an open real
analytic Riemannian manifold having the  same dimension as $\mathcal{M}$.
Due to the compactness of the manifold, the extended function $\hat{U}$ satisfies
\begin{align}
\triangle_g \hat{U}+\partial_t^2 \hat{U} +\langle B, \  \nabla \hat{U}\rangle +A \cdot \hat{U}=0 \quad \mbox{in}  \ \widehat{\mathcal{M}}_\rho\times (-\rho, \ -\rho),
\label{remind}
\end{align}
where $\widehat{\mathcal{M}}_\rho= \{ x\in \mathcal{M}_1| dist\{x, \ \mathcal{M}\}\leq \rho\}$.
From the  uniqueness of the analytic continuation, we also have that
\begin{align}
\triangle_g {U} +\langle B, \  \nabla {U}\rangle +A \cdot {U}+\lambda U= 0 \quad &\mbox{in} \ \widehat{\mathcal{M}}_\rho,
\label{fuzel}
\end{align}

It follows from (\ref{agree}) and the definition of $\hat{U}$ that
\begin{align}
\|U\|_{L^\infty(\mathbb B_{\rho})}\leq e^{C\sqrt{\lambda}} \|U\|_{L^\infty(\mathbb B^+_2)}.
\end{align}
See details of the similar arguments in e.g. \cite{LZ}.
Iterating the doubling inequality (\ref{LLL}) in the half balls finite number of steps, we can show that
\begin{align}
\|U\|_{L^\infty(\mathbb B_{\rho})}&\leq  e^{C\sqrt{\lambda}} \|U\|_{L^\infty(\mathbb B^+_\frac{\rho}{2})} \nonumber \\
&\leq e^{C\sqrt{\lambda}} \|U\|_{L^\infty(\mathbb B_\frac{\rho}{2})}
\end{align}
By rescaling arguments, it follows that
\begin{align}
\|U\|_{L^\infty(\mathbb B_{2r})}\leq  e^{C\sqrt{\lambda}} \|U\|_{L^\infty(\mathbb B_r)}
\label{correct}
\end{align}
for any $r\leq \frac{\rho}{2}$ with $\mathbb B_{2r}\subset \widehat{\mathcal{M}}_\rho$ and $C$ depending only on $\mathcal{M}$.

Next we need to extend $U(x)$ locally as a holomorphic function  in
$\mathbb{C}^n$. Applying elliptic estimates for $\hat{U}$ in (\ref{remind}) in a ball $\mathbb B_{r}(p)\times (-r, r)
\subset \widehat{\mathcal{M}}_\rho\times (-\rho, \rho)$ with $p\in \partial\mathcal{M}$ gives that
\begin{align}
|\frac{ D^{{\alpha}} \hat{U}(p, 0)}{{\alpha} !}|\leq
C^{|\alpha|}_1 r^{-|{\alpha}| } \|\hat{U}\|_{L^\infty},
\end{align}
where ${\alpha}$ is a multi-index taken with respect to $x$ and $C_1>1$ depends on $\mathcal{M}$.
By translation, we may consider the point $p$ as the origin. Recall the definition of $\hat{U}$ in (\ref{udefi}). We derive that
\begin{align}
|\frac{ D^{{\alpha}} {U}(0)}{{\alpha} !}|\leq
C^{|\alpha|}_1 r^{-|{\alpha}| } e^{C\sqrt{\lambda}} \|{U}\|_{L^\infty(\mathbb B_r(0)}.
\end{align}
We sum up a geometric series to extend $U(x)$ to be a
holomorphic function $U(z)$ with $z\in\mathbb{C}^n$. Thus, we have
\begin{align}
\sup_{|z|\leq \frac{r}{2C_1}}|U(z)|\leq C_2  e^{C\sqrt{\lambda}} \sup_{|x|\leq
r}|U(x)|
\end{align}
with $C_2>1$.
By iteration of the doubling inequality (\ref{correct}) finitely many times and the rescaling arguments, we derive that
\begin{align}
\sup_{|z|\leq 2r}|\nabla e_\lambda(z)|\leq e^{C_5\sqrt{\lambda}} \sup_{|x|\leq
r}|\nabla e_\lambda(x)| \label{dara}
\end{align}
for $0<r<\rho_0<\frac{\rho}{2}$, where $\rho_0$ and $C_5$ depend on $\mathcal{M}$ and $\mathbb B_r(0)\subset \widehat{\mathcal{M}}$.

Thanks to the doubling inequalities (\ref{dara}), we are ready to prove the upper bounds of critical sets of Dirichlet eigenfunctions as the arguments in Theorem \ref{th1}.
\begin{proof}[Proof of Theorem \ref{th2}] We first prove the critical sets in a
neighborhood of the boundary $\partial{\mathcal{M}}$. Let $G(x)=|\nabla e_\lambda(x)|^2$. We study the nodal sets of $G(x)$, which are just critical sets of $e_\lambda(x)$.
By rescaling and
translation, we can argue on scales of order one. Let $p\in \mathbb
B_{1/4}$ be the point where the maximum of $|\nabla e_\lambda|$ in $\mathbb
B_{1/4}$ is attained. For each fixed direction $\omega \in S^{n-1}$, set
${G}_{\omega}(z)=G(p+z\omega)$ in $z\in \mathcal{B}_1\subset\mathbb{C}$. Denote $N(\omega)=\sharp\{z\in\mathcal{B}_{1/2}(0)\subset
\mathbb{C}| G_\omega(z)=0\}$. Thanks to the doubling property (\ref{dara}) and the
Lemma \ref{wwhy}, we derive that
\begin{align}
\sharp\{ x \in \mathbb B_{1/2}(p) &|x-p \ \mbox{is parallel to} \
\omega \ \mbox{and} \ |\nabla e_\lambda(x)|=0\} \nonumber\\&\leq
\sharp\{z\in\mathcal{B}_{1/2}\subset
\mathbb{C}| G_\omega(z)=0\} \nonumber\\
&=N(\omega)\nonumber\\
&\leq C\sqrt{\lambda}.
\end{align}
 It follows from the integral geometry estimates that
\begin{align}
H^{n-1}(\{ x \in \mathbb B_{1/2}(p)|\ |\nabla e_\lambda(x)|=0\}) &\leq
c(n)\int_{S^{n-1}} N(\omega)\,d \omega \nonumber\\
&\leq \int_{S^{n-1}}C\sqrt{\lambda}\, d\omega \nonumber\\ &=C\sqrt{\lambda}.
\end{align}
Therefore, we arrive at
\begin{eqnarray}
H^{n-1}(\{ x \in \mathbb B_{1/4}| \ |\nabla e_\lambda(x)|=0\}) \leq C\sqrt{\lambda}.
\label{silin}
\end{eqnarray}
Since $\partial\mathcal{M}$ is compact, we can choose finitely many balls centered at $\partial\mathcal{M}$ so that those balls cover $\widehat{\mathcal{M}}_{\frac{\rho_0}{4}}$. From (\ref{silin}),
 we arrive at
\begin{equation}H^{n-1}(\{x\in
\widehat{\mathcal{M}}_{\frac{\rho_0}{4}}| \ |\nabla e_\lambda(x)|=0\})\leq C\sqrt{\lambda}.
\label{last12}
\end{equation}

Next we deal with the measure of nodal sets in
$\mathcal{M}\backslash \widehat{\mathcal{M}}_{\frac{\rho_0}{4}}$. We have
obtained the doubling inequality (\ref{III2}) in the interior of the manifold $\mathcal{M}$. We may choose $R\leq \frac{\rho_0}{16}$ in the interior doubling inequality (\ref{III2}).
 We also extend
$U(x)$ locally as a holomorphic function in $\mathbb {C}^n$.
Note that $U(x)$ satisfies (\ref{system1}) in $\mathcal{M}\backslash
\widehat{\mathcal{M}}_{\frac{\rho_0}{4}}$. We use the lifting argument as (\ref{udefi}) to get rid of $\lambda$. Thus, we have
\begin{align}
\triangle_g \hat{U}+\partial_t^2 \hat{U} +\langle B, \  \nabla \hat{U}\rangle +A \cdot \hat{U}=0 \quad \mbox{in}  \ \mathcal{M}\backslash
\widehat{\mathcal{M}}_{\frac{\rho_0}{4}}\times (-\infty, \ -\infty).
\end{align}
Applying elliptic estimates in a
small ball $\mathbb B_{r}(p)\times (-r, r)$, we have
\begin{equation}
|\frac{ D^\alpha \hat{U}(p, 0)}{\alpha !}|\leq
C^{|\alpha|}_3 r^{-|\alpha|} \|\hat{U}\|_{L^\infty},
\end{equation}
where $C_3>1$ depends only on $\mathcal{M}$ and $D^\alpha$ is taken with respect to $x$.
Let us consider the point $p$ as the origin as well. The definition of $\hat{U}$ further implies that
\begin{equation}
|\frac{ D^\alpha {U}( 0)}{\alpha !}|\leq
C^{|\alpha|}_3 r^{-|\alpha|}e^{C\sqrt{\lambda}} \|{U}\|_{L^\infty}.
\end{equation}
 By summing up a geometric series, we can extend $U(x)$ to be a
holomorphic function $U(z)$ with $z\in\mathbb{C}^n$ to have
\begin{align}
\sup_{|z|\leq \frac{r}{2C_3}}|U(z)|\leq C_4 e^{C\sqrt{\lambda}}\sup_{|x|\leq
r}|U(x)|
\end{align}
with $C_4>1$.
By finite steps of iterations of (\ref{III2}) and  rescaling arguments, we further derive that
\begin{equation}
\sup_{|z|\leq 2r}|U(z)|\leq e^{C \sqrt{\lambda}} \sup_{|x|\leq
r}|U(x)| \label{dara1}
\end{equation}
holds for $0<r<\frac{\rho_0}{16}$ with ${\rho_0}$ depending on $\mathcal{M}$. We make use of the Lemma \ref{wwhy} and the inequality
(\ref{dara1}) to obtain the upper bounds as the previous arguments  in the neighborhood of the
boundary.  By rescaling and translation, we can argue on scales
of order one. Let $p\in \mathbb B_{1/4}$ be the point where the
maximum of $|\nabla e_\lambda|$ in $\mathbb B_{1/4}$ is arrived. Recall that $G(x)=|\nabla e_\lambda(x)|^2.$ For each
direction $\omega \in S^{n-1}$, set $G_{\omega}(z)=
G(p+z\omega)$ in $z\in \mathcal{B}_1\subset\mathbb{C}$.
Thanks to the doubling property (\ref{dara1}) and the complex growth Lemma \ref{wwhy}, we derive that
\begin{align}
\sharp\{ x \in \mathbb B_{1/2}(p) &| x-p \ \mbox{is parallel to} \
\omega \ \mbox{and} \ |\nabla e_\lambda(x)|=0\} \nonumber\\&\leq
\sharp\{z\in\mathcal{B}_{1/2}\subset
\mathcal{C}| G_\omega(z)=0\} \nonumber\\
 &\leq C \sqrt{\lambda}.
\end{align}
From the integral geometry estimates, we have
\begin{align}
H^{n-1}(\{ x \in \mathbb B_{1/2}(p)|  \ |\nabla e_\lambda(x)|=0\}) &\leq
c(n)\int_{S^{n-1}} N(\omega)\,d \omega \nonumber\\
&\leq \int_{S^{n-1}}C\sqrt{\lambda}\, d\omega\nonumber\\ &=C\sqrt{\lambda}.
\end{align}
Thus, we deduce that
\begin{eqnarray}
H^{n-1}(\{ x \in \mathbb B_{1/4}| \ |\nabla e_\lambda(x)|=0\} )\leq C\sqrt{\lambda}.
\end{eqnarray}
 Covering the compact
manifold $\mathcal{M}\backslash\widehat{\mathcal{M}}_{\frac{\rho_0}{4}}$ with finite number of coordinate charts gives that
\begin{equation} H^{n-1}(\{x\in
\mathcal{M}\backslash\widehat{\mathcal{M}}_{\frac{\rho_0}{4}}| \ |\nabla e_\lambda(x)|=0\})\leq
C\sqrt{\lambda}. \label{last2}
\end{equation}
Together with the upper bounds in (\ref{last12}) and (\ref{last2}), we arrive at the
conclusion in Theorem \ref{th2}.

\end{proof}

\section{Proof of the Carleman estimates}
This section is devoted to the quantitative Carleman estimates in Proposition \ref{pro2}.
 \begin{proof}[Proof of Proposition \ref{pro2}]
We introduce the polar geodesic coordinates $ (r, \theta)$
in the half ball $\mathbb B^+_{r}$. Following the Einstein notation, we write
Laplace-Beltrami operator as
\begin{align*} r^2 \triangle =r^2 \partial^2_r + r^2\big( \partial_r \ln\sqrt{b}+\frac{n-1}{r}\big) \partial_r v+
\frac{1}{\sqrt{b}}
\partial_i\big(\sqrt{b}b^{ij}\partial_j \big),
\end{align*}
where $\partial_i =\frac{\partial}{\partial \theta_i}$, $
b_{ij}(r, \theta)$ is a metric on the geodesic sphere $S^{n-1}$, $b^{ij}$ is the inverse of $b_{ij}$,
$b=\det(b_{ij})$ and $\theta_{n-1}=\frac{x_n}{|x|}$. One can check that, for $r$ small
enough,
\begin{equation}
\left \{ \begin{array}{lll}
|\partial_r b^{ij}|\leq C| b^{ij}| \ \ \mbox{in term of
tensors},\medskip \\
|\partial_r b|\leq C, \label{gamma} \medskip\\
C^{-1} \leq b \leq C,
\end{array}
\right.
\end{equation}
where $C$ depends on ${\mathcal{M}}$. We want to transform the half ball $\mathbb B^+_{r}$ into a half cylinder.
Let $ r= e^t$. Then $\partial_r =e^{-t} \partial_t$. Hence the function $V(t, \theta_1, \cdots, \theta_{n-1})$ is supported in
$(-\infty, \  T_0]\times S^{n-1}_+$. Notice that $T_0$ is negative with large enough $|T_0|$ since $r$ is small. Under this new
coordinate, we can write
\begin{align} e^{2t} \triangle = \partial^2_t + ( n-2+ \partial_t\ln\sqrt{b}) \partial_t +
\frac{1}{\sqrt{b}}
\partial_i\big(\sqrt{b}b^{ij}\partial_j \big).
\label{new}
\end{align}
Furthermore, the condition (\ref{gamma})
turns into
\begin{equation}
\left \{ \begin{array}{lll}
|\partial_t b^{ij}|\leq C e^t |b^{ij}| \ \ \mbox{in term of
tensors},\medskip \\
|\partial_t b|\leq C e^t, \label{gamma1} \medskip\\
C^{-1} \leq b \leq C.
\end{array}
\right.
\end{equation}
The boundary conditions in (\ref{model}) in the new coordinates  change into
\begin{align}
V_1(t, \theta_1, \cdots, \theta_{n-2}, 0)=0, \cdots, V_{n-1}(t, \theta_1, \cdots, \theta_{n-2}, 0)=0,  \nonumber \\
\nu\cdot \nabla_{S^{n-1}} V_n= e^t h(t, \theta) V_n            \quad   \quad \mbox{if} \ \theta_{n-1}=0,
\end{align}
where $\nu$ is the unit outer normal on $S^{n-1}_+\cap \{x|x_n=0\}$. Let
\begin{align}
V= e^{-\beta \psi} W.
\end{align}
We introduce the conjugate operator,
\begin{align}\mathcal {L}_\beta (W)&=r^2  e^{\beta\psi(x)} \triangle (e^{-\beta \psi(x)}W) +  r^2 E(x) W \nonumber \\
&=e^{2t} e^{-\beta\phi(t)} \triangle (e^{\beta \phi(t)}W) +e^{2t} E(t, \theta)W.
\end{align}
By straightforward calculations,  it follows from (\ref{new}) that
\begin{align}\mathcal
{L}_\beta (W)&=\partial^2_t W+\big(2\beta \phi'+(n-2)+\partial_t \ln\sqrt{b}\big)\partial_t W
\nonumber \\ &+\big(\beta^2\phi'^2+\beta
\phi''+(n-2)\beta \phi'+\beta \partial_t \ln \sqrt{b}\phi'\big)
W+\triangle_{\theta} W+ e^{2t} E W ,
\end{align}
where $$\triangle_{\theta} W=\frac{1}{\sqrt{b}}\partial_i
(\sqrt{b} b^{ij}\partial_j W)$$
is the Laplace-Beltrami operator on $S^{n-1}$.
 We introduce the
following $L^2$ norm
$$ \|W\|^2_{\phi}=\int_{(-\infty, \ T_0]\times S^{n-1}_+} |W|^2\sqrt{b} \phi'^{-3} \, dt d\theta,           $$
where $d\theta$ is measure on $S^{n-1}$. We write the cylinder $S^{n-1}_+\times (-\infty, \ T_0]$ as $N=[0, 2\pi]\times [-\frac{\pi}{2}, \frac{\pi}{2}]\times \cdots \times [0, \frac{\pi}{2}]\times (-\infty, \ T_0]$. Then $\partial S^{n-1}_+\cap \{x|x_n=0\}\times (-\infty, \ T_0]$ is denoted as $\partial N=[0, 2\pi]\times [-\frac{\pi}{2}, \frac{\pi}{2}]\times \cdots \{0\}\times (-\infty, \ T_0]$.

To obtain the Carleman estimates in the proposition, we aim to find a lower estimate for
$\| {L}_\beta (W)\|_{\phi}$ using some elementary algebra and integration by parts arguments. Note that $t$ is close to $-\infty$ in the later calculations.
By the triangle inequality, we easily have
$$ \|\mathcal{L}_\beta (W)\|_\phi\geq \mathcal{A}-\mathcal{B},              $$
where \begin{align} \mathcal{A}=\|
\partial^2_t W +2\beta \phi'\partial_tW+\beta^2 \phi'^2 W+ e^{2t}EW+\triangle_{\theta} W \|_{\phi}
 \end{align} and
 \begin{align}
\mathcal{B}=\|\beta \phi'' W+(n-2)\beta \phi' W+ \beta \partial_t \ln\sqrt{b} \phi'
W+(n-2)\partial_t W+\partial_t \ln\sqrt{b}\partial_t
W\|_{\phi}.
\end{align}
Later on, we can show that $\mathcal{B}$ can be controlled by $\mathcal{A}$. Thus, our main goal is to find a lower bound for  $\mathcal{A}$. Write
\begin{align}
 \mathcal{A}^2= \mathcal{A}_1+ \mathcal{A}_2+ \mathcal{A}_3,
 \label{cpp}
\end{align}
where
\begin{align}
\mathcal{A}_1=\|\partial^2_t W +( \beta^2 \phi'^2+ e^{2t}E)W+ \triangle_{\theta} W \|^2_{\phi},
\end{align}
\begin{align}
\mathcal{A}_2=\|2\beta \phi'\partial_t W\|^2_{\phi},
\label{AAA2}
\end{align}
\begin{align}
\mathcal{A}_3= 2\langle 2\beta \phi'\partial_tW, \ \, \partial^2_t W +\beta^2 \phi'^2W+ e^{2t}EW+\triangle_{\theta} W \rangle_{\phi}.
\end{align}
General speaking, we compute the Carleman estimates by writing $ \mathcal{A}$ as symmetric and antisymmetric parts. For these quantitative Carleman estimates, it save more computations if we suppress the terms with less contribution on $\beta$. Let us first compute the contribution from $\mathcal{A}_3$.  We decompose the inner product
$\mathcal{A}_3$ as follows
\begin{align*}
\mathcal{A}_3=K_1+K_2+K_3,
\end{align*}
where each integration $K_i$ is
\begin{align*}
K_1=4\beta \int_N \phi' \partial_t W \partial^2_t W  \phi'^{-3}\sqrt{b} \,dt  d \theta,
\end{align*}
\begin{align*}
K_2=4\beta \int_N \phi'\partial_t W  \partial_i(\sqrt{b} b^{ij}\partial_j W)  \phi'^{-3} \, dt  d \theta,
\end{align*}
\begin{align*}
K_3=4  \int_N (\beta^3 \phi'^{2}+\beta e^{2t} E) \partial_t W W  \phi'^{-2} \sqrt{b} \, dt d\theta.
\end{align*}

For $K_1$, applying the integration by parts with respect to $t$ gives that
\begin{align}
K_1=4\beta \int_N \phi''  |\partial_t W|^2\phi'^{-3}  \sqrt{b} dt\ d\theta -2\beta \int_N  \phi' |\partial_t W|^2 \partial_t \ln \sqrt{b}\phi'^{-3} \sqrt{b} \,dt  d \theta.
\end{align}
It follows from (\ref{gamma1}) that we have $|\partial_t \ln  \sqrt{b}|\leq C e^t$. Considering the definition of $\phi$, since $T_0$ close to $-\infty$, the term
$|\partial_t \ln  \sqrt{b}|$ is controlled by $|\phi''|=\frac{2}{t^2}$ and $\phi'\approx 1$. Thus, we have
\begin{align}
K_1\geq -C\beta \int_N |\phi''|  |\partial_t W|^2 \sqrt{b}\phi'^{-3} \, dt d\theta.
\label{mamabi}
\end{align}
To estimate $K_2$, we first integrate with respect to $\partial_i$ with the consideration that $W\in C^{\infty}_{0}(\mathbb B^+_{r_0})$,
\begin{align}
K_2&=-4\beta \int_N  \phi' \partial_t \partial_i W b^{ij} \partial_j W  \phi'^{-3} \sqrt{b} \, dt d\theta
+4\beta \int_{\partial N}
\phi' \partial_t W \partial_j W b^{j(n-1) }  \phi'^{-3} \sqrt{b} \, dt d\theta. \nonumber \\
&=J_1+J_2,
\label{kk2}
\end{align}
where
\begin{align*}
J_1=-4\beta \int_N  \phi' \partial_t \partial_i W b^{ij} \partial_j W  \phi'^{-3} \sqrt{b} \, dt d\theta
\end{align*}
and the boundary term
\begin{align*}
J_2=4\beta \int_{\partial N}
\phi' \partial_t W \partial_j W b^{j(n-1) }  \phi'^{-3} \sqrt{b} \, dt d\theta.
\end{align*}

From the boundary condition on $V$, we obtain that
\begin{align}
\partial_t W_j= \frac{\partial W_j}{\partial \theta_k}=0, \ \mbox{for}  \ 1\leq j\leq n-1, \ 1\leq k\leq n-2&  \ \ \mbox{for} \ \theta_{n-1}=0, \nonumber \\
\nu \cdot \nabla_{S^{n-1}} W_n=h e^t W_n& \ \mbox{for} \ \ \theta_{n-1}=0.
\label{boundary}
\end{align}
Thus, the boundary term in (\ref{kk2}) can be converted into
\begin{align}
J_2&=4\beta \int_{ \partial N}\phi' \partial_t W_n  W_n h e^t  \phi'^{-3} \sqrt{b} \, dt d\theta\nonumber \\
&=-2\beta \int_{\partial N}| W_n|^2  \partial_t( he^t  \phi'^{-2} \sqrt{b}) \, dt d\theta,
\end{align}
where we have used the integration by parts in the last identity with respect to $t$. Then
\begin{align}
|J_2|\leq \beta (\|h\|_{C^1}+1)\int_{\partial N}| W_n|^2    e^t  \phi'^{-3} \sqrt{b} \, dt d\theta.
\end{align}

To deal with the contribution from the boundary, the following trace lemma is established in \cite{R},
\begin{align}
\| u\|^2_{S^{n-2}}\leq C\big ( \tau  \| u\|^2_{S^{n-1}_+}+ \tau^{-1}  \| \nabla u\|^2_{S^{n-1}_+}\big)
\label{trace}
\end{align}
for any $\tau\geq 1$ and for $u\in C^\infty(S^{n-1}_+) $. The inequality (\ref{trace}) is obtained by the trace lemma in $\mathbb R^{n+1}$ and rescaling arguments.  From the trace inequality (\ref{trace}), we have
\begin{align}
\int_{\partial N}| W_n|^2    e^t  \phi'^{-3} \sqrt{b} \, dt d\theta &\leq C \tau \int^{T_0}_{-\infty} \int_{ S^{n-1}_+}| W_n|^2    e^t  \phi'^{-3} \sqrt{b} \, dt d\theta \nonumber \\ &+ C \tau^{-1}\int^{T_0}_{-\infty} \int_{ S^{n-1}_+}|\nabla_\theta W_n|^2    e^t  \phi'^{-3} \sqrt{b} \, dt d\theta.
\end{align}
Thus, we choose
\begin{align}
\tau=\beta\geq  C(1+\sqrt{\|E\|}_{C^1}+\|h\|_{C^1})
\end{align}
to have
\begin{align}
|J_2|\leq \beta^3 \int_{ N}| W_n|^2    e^t  \phi'^{-3} \sqrt{b} \, dt d\theta + \beta \int_{ N}|\nabla_\theta W_n|^2    e^t  \phi'^{-3} \sqrt{b} \, dt d\theta.
\end{align}

Integrating by parts with respect to $t$ for $J_1$ gives that
\begin{align}
J_1=&-4\beta \int_N  \phi''  \partial_i W b^{ij} \partial_j W  \phi'^{-3} \sqrt{b} \, dt d\theta
+ 2\beta \int_N  \phi'  \partial_t \ln \sqrt{b} \partial_i W b^{ij} \partial_j W  \phi'^{-3} \sqrt{b} \, dt d\theta \nonumber \\
&+2\beta \int_N  \phi'  \partial_t b^{ij} \partial_i W  \partial_j W  \phi'^{-3} \sqrt{b} \, dt d\theta.
\end{align}
Denote that
\begin{align*}
|\nabla_\theta W|^2= b^{ij} \partial_i W  \partial_j W.
\end{align*}
From the fact that $-\phi''$ is nonnegative, $\beta$ is large and the assumption for $b^{ij}$, we arrive at
\begin{align}
J_1\geq 3 \beta \int_N  |\phi'' | |\nabla_\theta W|^2 \phi'^{-3} \sqrt{b} \, dt d\theta.
\end{align}
Combining the estimates for $J_1$ and $J_2$ and taking into account that $|\phi''|\gg e^t$ for $|T_0|$ large enough yield that
\begin{align}
K_2\geq 2\beta \int_N  |\phi'' | |\nabla_\theta W|^2 \phi'^{-3} \sqrt{b} \, dt d\theta- C\beta^3 \int_N  |\phi'' | | W|^2 \phi'^{-3} \sqrt{b} \, dt d\theta.
\label{mamabi1}
\end{align}

We carry out the similar computations for $K_3$,
\begin{align}
K_3=&-  2\beta^3 \int_N |W|^2  \partial_t \ln \sqrt{b}    \sqrt{b} \, dt d\theta \nonumber \\& - \int_N ( 4 \phi'-4\phi^{''}+ 2\phi'
\partial_t \ln \sqrt{b} )\beta e^{2t} E |W|^2 \phi'^{-3} \sqrt{b} \, dt d\theta \nonumber \\
&- 2\beta \int_N \phi'  e^{2t} \partial_t E   |W|^2  \phi'^{-3} \sqrt{b} \, dt d\theta.
\end{align}
Since we have assumed that
\begin{align*}
\beta>C(1+\sqrt{\|E\|}_{C^1}+\|h\|_{C^1}),
\end{align*}
from the fact that $|T_0|$ is large enough, (\ref{gamma1}) and $\phi'\approx 1$, we obtain that
\begin{align}
K_3\geq - C \beta^3 \int_N |W|^2    e^t   \phi'^{-3} \sqrt{b} \, dt d\theta.
\label{mamabi2}
\end{align}
Combining the estimates (\ref{mamabi}) for $K_1$, (\ref{mamabi1}) for $K_2$ and (\ref{mamabi2}) for $K_3$ yields that
\begin{align}
\mathcal{A}_3& \geq 2\beta \int_N |\phi''|  |\nabla_\theta W|^2 \phi'^{-3} \sqrt{b} \, dt d\theta -C\beta^3
\int_N e^t | W|^2 \phi'^{-3} \sqrt{b} \, dt d\theta \nonumber \\
&-C\beta \int_N |\phi''| | \partial_t W|^2 \phi'^{-3} \sqrt{b} \, dt d\theta.
\label{A33}
\end{align}

We need to obtain a stronger $L^2$ norm of $W$. To this end, we consider $\mathcal{A}_1$. Choose some small positive constant $\delta$ which is to be determined later. Since
$|\phi''|\leq 1$ and $\beta\geq 1$,  it follows that
\begin{align}
\mathcal{A}_1\geq \frac{\delta}{\beta}\hat{\mathcal{A}}_1,
\label{mm1}
\end{align}
where $\hat{\mathcal{A}}_1$ is given by
\begin{align}
\hat{\mathcal{A}}_1=\||\phi''|^{\frac{1}{2}}  \big(\partial^2_t W +( \beta^2 \phi'^2+ e^{2t}E)W+ \triangle_{\theta} W \big)\|^2_{\phi}.
\end{align}
Furthermore, we decompose $\hat{\mathcal{A}}_1$ as
\begin{align}
\hat{\mathcal{A}}_1= \mathcal{K}_1+\mathcal{K}_2+\mathcal{K}_3
\label{mm2}
\end{align}
where
\begin{align*}\mathcal{K}_1=\| |\phi''|^{\frac{1}{2}} \big(\partial^2_t W+\triangle_{\theta} W \big)\|_\phi^2 \end{align*} and
\begin{align*}\mathcal{K}_2=\|  |\phi''|^{\frac{1}{2}} \big(\beta^2 \phi'^2+
e^{2t}E \big) W\|_\phi^2 \end{align*}
and
\begin{align}\mathcal{K}_3=2\langle |\phi^{''}|(\partial^2_t W+\triangle_{\theta} W), \ \big(\beta^2 \phi'^2+
e^{2t} E \big) W \rangle_\phi.
\label{koko} \end{align}
We further split $\mathcal{K}_3$ into
$\mathcal{K}_3=H_1+H_2$,
where
\begin{align*}
H_1 &= 2 \int_N  \partial^2_t W  \big(\beta^2 \phi'^2+
e^{2t} E \big) W |\phi^{''}| \phi'^{-3} \sqrt{b}\, dt d\theta
\end{align*}
and
\begin{align*}
H_2= 2 \int_N |\phi^{''} |\triangle_{\theta} W \big(\beta^2 \phi'^2+
e^{2t} E \big) W \phi'^{-3} \sqrt{b} \, dt d\theta.
\end{align*}

Integration by parts with respect to $t$ gives that
\begin{align}
H_1
= & \int_N \phi^{''}   \big(\beta^2 \phi'^2+
e^{2t} E \big)   |\partial_t W|^2 \phi'^{-3}   \sqrt{b} \, dt d\theta\nonumber \\
 &+ 2 \int_N \partial_t \big[(\beta^2 \phi'^2+
e^{2t} E ) \phi^{''} \phi'^{-3} \sqrt{b}\big] W  \partial_t W \, dt d\theta.
\end{align}
By Cauchy-Schwartz inequality and the assumption of $\beta$, we obtain that
\begin{align}
H_1\geq -C\beta^2 \int_N |\phi^{''}|  (  |\partial_t W|^2 +|W|^2) \phi'^{-3}   \sqrt{b} \, dt d\theta.
\label{HH1}
\end{align}

We perform the integration by part arguments with respect to $\partial_i$ gives that
\begin{align}
H_2=&2 \int_N \phi^{''}  |\nabla_{\theta} W|^2  \big(\beta^2 \phi'^2+
e^{2t} E \big)  \phi'^{-3}\sqrt{b}\,dt d\theta  \nonumber \\
&+2 \int_N \phi^{''}e^{2t}  \partial_i W b^{ij}\partial_j EW \phi'^{-3} \sqrt{b}\,dt d\theta \nonumber \\
&-2 \int_{\partial N} \phi^{''}e^{t} \big(\beta^2 \phi'^2+
e^{2t} E \big)   h W^2_n\phi'^{-3} \sqrt{b}\,dt d\theta,
\end{align}
where we have used the boundary conditions (\ref{boundary}). From Cauchy-Schwartz inequality and the assumption of $\beta$, it holds that
\begin{align*}
|\partial_ jE b^{ij} \partial_i W W|\leq C\beta^2 (|\nabla_\theta W|^2+|W|^2).
\end{align*}
Therefore, we further obtain that
\begin{align}
H_2\geq &-C\beta^2 \int_N |\phi^{''}| (|\nabla_{\theta} W|^2 +|W|^2) \phi'^{-3}\sqrt{b}\,dt d\theta  \nonumber \\
&-C\|h\|_{C^1} \beta^2 \int_{\partial N} |\phi^{''}| e^{2t} W^2_n\phi'^{-3} \sqrt{b}\,dt d\theta.
\label{HH2}
\end{align}
With the aid of the trace inequality (\ref{trace}) by choosing $\tau=\beta$ and the assumption of $\beta$, we derive that
\begin{align}
\|h\|_{C^1} \beta^2 \int_{\partial N} |\phi^{''}| e^{2t} W^2_n\phi'^{-3} \sqrt{b}\,dt d\theta &\leq C\beta^4   \int_{ N} |\phi^{''}| e^{2t} W^2_n\phi'^{-3} \sqrt{b}\,dt d\theta \nonumber \\&+C\beta^2 \int_{ N} |\phi^{''}| e^{2t}|\nabla_\theta W_n|^2\phi'^{-3} \sqrt{b}\,dt d\theta.
\label{HH3}
\end{align}
Thus, we learn from (\ref{koko}), (\ref{HH1}), (\ref{HH2}) and (\ref{HH3}) that
\begin{align}
\mathcal{K}_3&\geq -C \beta^2 \int_N |\phi^{''}| \big(|\partial_t W|^2+ |\nabla_{\theta} W|^2+ |W|^2  \big) \phi'^{-3}\sqrt{b} \,dt d\theta \nonumber \\
&-C\beta^4  \int_{ N}| \phi^{''}| e^{2t}   |W|^2 \phi'^{-3} \sqrt{b}\,dt d\theta.
\label{mm3}
\end{align}
Since $\phi'$ is close to $1$ and $e^{2t}$ is sufficiently small as $|T_0|$ is large enough, it follows that
\begin{align}
\mathcal{K}_2\geq C\beta^4 \int_N |\phi''|  |W|^2   \phi'^{-3}\sqrt{b} \,dt d\theta.
\label{K22}
\end{align}

Note that $\mathcal{K}_1\geq 0$. If follows from (\ref{mm1}), (\ref{mm2}), (\ref{mm3}) and (\ref{K22}) that
\begin{align}
\mathcal{A}_1
\geq& -C\beta \delta \int_N |\phi^{''}| (|\partial_t W|^2+ |\nabla_{\theta} W|^2 )\phi'^{-3} \sqrt{b} \,dt d\theta \nonumber \\
&+C\beta^3 \delta \int_N |\phi^{''} |W|^2  \phi'^{-3} \sqrt{b}\,dt d\theta.
\label{A11}
\end{align}

Combining the estimates from (\ref{cpp}), (\ref{AAA2}), (\ref{A33}) and (\ref{A11}), we derive that
\begin{align}
\mathcal{A}^2&\geq C\beta^3 \delta \int_N |\phi^{''}| | W|^2 \phi'^{-3}\sqrt{b}\,dt d\theta +4\beta^2 \int_N |\phi^{'}| |\partial_t W|^2 \phi'^{-3} \sqrt{b}\,dt d\theta   \nonumber \\
 &+2\beta \int_N |\phi^{''}| |\nabla_{\theta} W|^2 \phi'^{-3}\sqrt{b}\,dt d\theta
- C\beta^3 \int_N e^t | W|^2 \phi'^{-3}\sqrt{b}\,dt d\theta \nonumber \\
& -C \beta \delta\int_N |\phi^{''}|\big( |\partial_t W|^2+ |\nabla_{\theta} W|^2  )  \phi'^{-3} \sqrt{b}\,dt d\theta\nonumber \\
&  -C\beta \int_N |\phi^{''}| |\partial_t W|^2  \phi'^{-3} \sqrt{b} \,dt d\theta.
\end{align}
Since $\delta$  can be chose to be small and $e^t$ is small enough compared with $|\phi^{''}| $, it follows that
\begin{align}
\mathcal{A}^2&\geq C\beta^3  \int_N |\phi^{''}| | W|^2 \phi'^{-3}\sqrt{b}\,dt d\theta
+C\beta^2 \int_N |\phi^{'}| |\partial_t W|^2\phi'^{-3}\sqrt{b}  \,dt d\theta
  \nonumber \\
 &+C\beta \int_N |\phi^{''}| |\nabla_{\theta} W|^2 \phi'^{-3}\sqrt{b}\,dt d\theta.
 \label{lower}
\end{align}
It is easy to check that we can absorb $\mathcal{B}$ to $\mathcal{A}$ if $|T_0|$ is large enough and $\beta$ is chosed to be large.  Thus, we obtain that
\begin{align}
\|\mathcal{L}_\beta (W)\|^2_\phi  &\geq C\beta^3  \int_N |\phi^{''}| | W|^2\phi'^{-3}\sqrt{b} \,dt d\theta
+C\beta^2 \int_N |\phi^{'}| |\partial_t W|^2 \phi'^{-3}\sqrt{b} \,dt d\theta
  \nonumber \\
 &+C\beta \int_N |\phi^{''}| |\nabla_{\theta} W|^2 \phi'^{-3}\sqrt{b}\,dt d\theta.
\end{align}

Since $|\phi''|$ is much smaller compared with $|\phi'|$, we have that
\begin{align}
\|\mathcal{L}_\beta (W)\|^2_\phi  &\geq C\beta^3  \int_N |\phi^{''}| | W|^2 \phi'^{-3}\sqrt{b}\,dt d\theta
+\hat{C}\beta \int |\phi^{''}| |\partial_t W|^2 \phi'^{-3} \sqrt{b}\,dt d\theta
  \nonumber \\
 &+C\beta \int_N |\phi^{''}| |\nabla_{\theta} W|^2\phi'^{-3}\sqrt{b} \,dt d\theta,
\end{align}
where $\hat{C}$ can be chosen arbitrarily smaller than $C$. Recall that the conjugate operator $V= e^{-\beta\psi} W$. It leads to
\begin{align}
\| e^{2t} e^{\beta \psi} (\triangle V+ EV)\|^2_\phi  &\geq C\beta^3  \| |\phi^{''}|^{\frac{1}{2}}e^{\beta \psi}   V\|_{\phi}^2
+\hat{C}\beta \| |\phi^{''}|^{\frac{1}{2}}e^{\beta \psi}  | \partial_t V|\|_{\phi}^2
  \nonumber \\
 &- \hat{C}\beta^3 \| |\phi^{''}|^{\frac{1}{2}}e^{\beta \psi}  V\|_{\phi}^2
 +C\beta \| |\phi^{''}|^{\frac{1}{2}} e^{\beta \psi} |\nabla_{\theta} V|\|_{\phi}^2.
\end{align}
Since $\hat{C}$ is chose to be smaller than $C$, we derive that
\begin{align}
\| e^{2t} e^{\beta \psi} (\triangle V+ EV)\|^2_\phi  &\geq C\beta^3  \| |\phi^{''}|^{\frac{1}{2}}e^{\beta \psi}  V\|_{\phi}^2
+\hat{C}\beta \| |\phi^{''}|^{\frac{1}{2}}e^{\beta \psi}  | \partial_t V|\|_{\phi}^2
  \nonumber \\
 &+C\beta \| |\phi^{''}|^{\frac{1}{2}} |\nabla_{\theta} V\|_{\phi}^2.
 \label{volume}
\end{align}
Note that the volume element in polar coordinates is $r^{n-1} \sqrt{b}  dr d\theta$ and $\phi'\approx 1$. It follows from (\ref{volume}) that
\begin{align}
\| r^2 e^{\beta \psi} (\triangle V+ EV) r^{-\frac{n}{2}}\|^2_{\mathbb B^+_{r_0}}  &\geq C\beta^3  \| (\log r)^{-1}e^{\beta \psi}   Vr^{-\frac{n}{2}}\|^2_{\mathbb B^+_{r_0}}  \nonumber \\
&+{C}\beta \|(\log r)^{-1}e^{\beta \psi}  | \nabla V|r^{-\frac{n}{2}} \|^2_{\mathbb B^+_{r_0}}.
\end{align}
By replacing $\beta$ by $\beta+\frac{n}{2}$, which only changes the lower bound of $\beta$ by a constant in the Carleman estimates, we arrive at the desired estimates (\ref{Carle}).

Next we prove  stronger Carleman estimates with some strong assumption on $V$.
Suppose that $\supp V\subset \{ x\in \mathbb B^+_{r_0}| r(x)\geq \rho\}$. Let $\hat{T}_0=\ln \rho$. The application of Cauchy-Schwarz inequality gives that
\begin{align}
\int_N \partial_t |W|^2 e^{-t} \sqrt{b} \, dtd\theta\leq 2 (\int_N |\partial_t W|^2 e^{-t} \sqrt{b} \, dtd\theta)^{\frac{1}{2}}
(\int_N | W|^2 e^{-t} \sqrt{b} \, dtd\theta)^{\frac{1}{2}}.
\label{cauchy}
\end{align}
For the left hand side of (\ref{cauchy}), applying the integration by parts shows that
\begin{align}
\int_N \partial_t |W|^2 e^{-t} \sqrt{b} \, dtd\theta=\int_N |W|^2 e^{-t} \sqrt{b} \, dtd\theta-\int_N |W|^2 e^{-t} \partial_t (\ln\sqrt{b})\sqrt{b} \, dtd\theta.
\end{align}
Since $|\partial_t\ln\sqrt{b}|\leq C e^t$ for $|\hat{T}_0|$ large enough, we obtain that
\begin{align}
\int_N \partial_t |W|^2 e^{-t} \sqrt{b} \, dtd\theta\geq C \int_N  |W|^2 e^{-t} \sqrt{b} \, dtd\theta.
\label{cauchy1}
\end{align}
Taking (\ref{cauchy}) and (\ref{cauchy1}) into consideration gives that
\begin{align}
e^{-\hat{T}_0}\int_N | \partial_t W|^2  \sqrt{b} \, dtd\theta &\geq \int_N | \partial_t W|^2 e^{-t} \sqrt{b} \, dtd\theta \nonumber \\
&\geq C\int_N |W|^2 e^{-t} \sqrt{b} \, dtd\theta.
\end{align}
Notice that $e^{-\hat{T}_0}=\rho^{-1}$. From (\ref{lower}), we deduce that
\begin{align}
\mathcal{A}^2\geq C\beta^2 \rho \int_N |W|^2 e^{-t}  \phi'^{-3}\sqrt{b} \, dtd\theta.
\label{latter}
\end{align}
Thus, together with (\ref{lower}) and (\ref{latter}), we arrive at
\begin{align}
\mathcal{A}^2&\geq C\beta^3  \int_N |\phi^{''}| | W|^2 \phi'^{-3}\sqrt{b}\,dt d\theta
+C\beta^2 \int_N |\phi^{'}| |\partial_t W|^2 \phi'^{-3} \sqrt{b}\,dt d\theta
  \nonumber \\
 &+C\beta \int_N |\phi^{''}| |\nabla_{\theta} W|^2 \phi'^{-3}\sqrt{b}\,dt d\theta +C\beta^2 \rho \int_N |W|^2 e^{-t}  \phi'^{-3}\sqrt{b} \, dtd\theta.
 \label{ahead}
\end{align}
Therefore, the Carleman estimates (\ref{Strongcarle}) follows from (\ref{ahead}) as the deduction of (\ref{Carle}) in the previous arguments.
\end{proof}


\begin{thebibliography}{CL}

\bibitem[ARRV]{ARRV} G. Alessandrini, L. Rondi, E. Rosset and S. Vessella,
The stability for the Cauchy problem for elliptic equations.
Inverse Problems,  25(2009), no. 12, 123004, 47 pp.

\bibitem[BLS]{BLS} L. Buhovsky, A. Logunov and M. Sodin, Eigenfunctions with infinitely many isolated critical points, To appear in International Mathematics Research Notices.

\bibitem[B]{B}L. Bakri, Critical sets of eigenfunctions of the Laplacian, Journal of Geometry and Physics, 62(2012), no.10, 2024-2037.

\bibitem[BL]{BL}K. Bellova and F.-H. Lin, Nodal sets of Steklov eigenfunctions,
Calc. Var. $\&$ PDE, 54(2015), 2239--2268.

\bibitem[CNV]{CNV} J. Cheeger, A. Naber and D. Valtorta, Critical sets of elliptic equations, Comm. Pure Appl. Math, 68(2015), no.2, 173-209.

\bibitem[DF]{DF}H. Donnelly and C. Fefferman, Nodal sets of eigenfunctions
on Riemannian manifolds, Invent. Math. 93(1988), 161-183.

\bibitem[DF1]{DF1} H. Donnelly and C. Fefferman, Nodal sets of eigenfunctions:
Riemannian manifolds with boundary, in: Analysis, Et Cetera,
Academic Press, Boston, MA, 1990, 251--262.

\bibitem[JN]{JN} D. Jakobson and N. Nadirashvili, Eigenfunctions with few critical points, J. Differential Geometry, 53(1999), 177-182.


\bibitem[Han]{Han} Q. Han,  Singular sets of solutions to elliptic equations, Indiana Univ. Math. J., 43(1994), no.3, 983-1002.

\bibitem[Han1]{Han1} Q. Han,  Singular sets of harmonic functions in R2 and their complexifications in $\mathbb C^2$,  Indiana Univ. Math. J., 53(2004), no.5, 1365-1380.

\bibitem[HL]{HL} Q. Han and F.-H. Lin, Nodal sets of solutions of Elliptic
Differential Equations, book in preparation (online at
http://www.nd.edu/qhan/nodal.pdf).

\bibitem[HHL]{HHL} Q. Han, R. Hardt and F.-H. Lin, Geometric measure of singular sets of elliptic equations, Comm. Pure Appl. Math., 51(1998), no.11-12, 1425-1443.

\bibitem[HHON]{HHON} R. Hardt, M. Hoffmann-Ostenhof, T. Hoffmann-Ostenhof and N. Nadirashvili, Critical sets of
solutions to elliptic equations, J. Differential Geom, 51(1999), no.2, 359-373.


\bibitem[Lin]{Lin}F.-H. Lin, Nodal sets of solutions of elliptic
equations of elliptic and parabolic equations, Comm. Pure Appl Math.
44(1991), 287-308.

\bibitem[LZ]{LZ} F.H. Lin and J. Zhu, Upper bounds of nodal sets for eigenfunctions of eigenvalue problems, arXiv:2005.04079.


\bibitem[Lo1]{Lo1} A. Logunov, Nodal sets of Laplace eigenfunctions: polynomial upper estimates of the Hausdorff measure,  Ann. of Math., 187(2018), 221--239.

\bibitem[Lo2]{Lo2} A. Logunov, Nodal sets of Laplace eigenfunctions: proof of Nadirashvili's conjecture and of the lower bound in Yau's conjecture, 187(2018), 241--262.

\bibitem[LM]{LM} A. Logunov and E. Malinnikova, Nodal sets of Laplace eigenfunctions: estimates of the Hausdorff measure in dimension two and three, 	50 years with Hardy spaces, 333-344, Oper. Theory Adv. Appl., 261, Birkh\"auser/Springer, Cham, 2018.


\bibitem[NV]{NV}A. Naber and D. Valtorta, Volume estimates on the critical sets of solutions to elliptic PDEs, Comm. Pure Appl. Math. 70(2017), no. 1835-1897.

\bibitem[Mo]{Mo} C. Morrey, Multiple integrals in the calculus of variations, Reprint of the 1966 edition,Springer-Verlag, Berlin, 2008, x+506 pp.


\bibitem[MN]{MN} C.B. Morrey and L. Nirenberg, On the analyticity of the solutions of linear elliptic systems of partial differential equations, 10(1957) 271-290.

\bibitem[R]{R} A. R\"uland,  quantitative unique continuation properties of fractional Schr\"odinger equations: doubling, vanishing order and nodal domain estimates, Trans. Amer. Math. Soc., 369(2017), no.4, 2311-2362.

\bibitem[TZ]{TZ}J. Toth and S. Zelditch,
Counting nodal lines which touch the boundary of an analytic domain,
J. Differential Geom., 81(2009), no.3, 649-686.

\bibitem[TZ1]{TZ1}J. Toth and S. Zelditch, Nodal intersections and Geometric Control, arXiv:1708.05754.


\bibitem[Y]{Y} S.-T. Yau, Problem section, Seminar on Differential Geometry, Annals of Mathematical Studies
102, Princeton, 1982, 669-706.

\bibitem[Z]{Z} S. Zelditch, Local and global analysis of eigenfunctions on Riemannian manifolds, in Handbook of geometric analysis, No. 1, vol. 7 of Adv. Lect. Math. (ALM), Int. Press, Somerville, MA, 2008, 545-658. 

\bibitem[Z1]{Z1} S. Zelditch, Measure of nodal sets of analytic
Steklov eigenfunctions, Math. Res. Lett., 22(2015), no.6, 1821--1842.

\bibitem[Zh]{Zh}J. Zhu, Doubling property and vanishing order of Steklov
eigenfunctions, Comm. Partial Differential Equations 40(2015), no.
8, 1498-1520.

\bibitem[Zh1]{Zh1}J. Zhu, Interior nodal sets of Steklov eigenfunctions on surfaces,
 Anal. PDE, 9(2016), no. 4, 859--880.

\bibitem[Zh2]{Zh2}J. Zhu, Geometry and interior nodal sets of Steklov eigenfunctions,   Calculus of
Variations and Partial Differential Equations,  59(2020), no. 5, 150.

\bibitem[Zh3]{Zh3} J. Zhu, Nodal sets of Robin and Neumann eigenfunctions,  arXiv:1810.12974.

\end{thebibliography}
\end{document}